\documentclass{amsproc}
\usepackage{hyperref}
\usepackage{amsmath}
\renewcommand{\enspace}{\;}
\catcode`\@=11
\def\@begintheorem#1#2[#3]{%
  \item[]{\normalfont 
  \hskip\labelsep
  \the\thm@headfont
  \thm@indent
  \@ifempty{#1}{\let\thmname\@gobble}{\let\thmname\@iden}%
  \@ifempty{#2}{\let\thmnumber\@gobble}{\let\thmnumber\@iden}%
  \@ifempty{#3}{\let\thmnote\@gobble}{\let\thmnote\@iden}%
  \thm@swap\swappedhead\thmhead{#1}{#2}{#3}%
  \the\thm@headpunct}%
  \@restorelabelsep
  \thmheadnl 
  \ignorespaces}
\catcode`\@=12

\DeclareMathAlphabet{\mathbbold}{U}{bbold}{m}{n}
\usepackage{mymacros}
\usepackage{graphicx}
\usepackage{color}
\newtheorem{theorem}{Theorem}[section]
\newtheorem{lemma}[theorem]{Lemma}
\newtheorem{corollary}[theorem]{Corollary}
\newtheorem{proposition}[theorem]{Proposition}
\theoremstyle{definition}
\newtheorem{definition}[theorem]{Definition}
\newtheorem{example}[theorem]{Example}
\theoremstyle{remark}
\newtheorem{remark}[theorem]{Remark}
\title{Max-plus convex sets and functions}
\author{Guy Cohen}
\address{Guy Cohen: 
Cermics-ENPC, 77455 Marne-La-Vall\'{e}e, cedex 2, France.}
\email{Guy.Cohen@mail.enpc.fr}
\author{St\'ephane Gaubert}
\address{St\'ephane Gaubert: 
INRIA-Rocquencourt, 78153 Le Chesnay cedex, France.}
\email{Stephane.Gaubert@inria.fr}
\author{Jean-Pierre Quadrat}
\address{Jean-Pierre Quadrat:
INRIA-Rocquencourt, 78153 Le Chesnay cedex, France.}
\email{Jean-Pierre.Quadrat@inria.fr}
\author{Ivan Singer}
\address{Ivan Singer: 
 Institute of Mathematics of the Romanian Academy, Bucharest, 
70700, Romania.}
\email{ivan.singer@imar.ro}

\thanks{This work was partially supported by the Erwin Schr\"odinger International Institute for Mathematical Physics (ESI) and the CERES program of the Romanian Ministry of Education and Research, contract no. 152/2001.}
\thanks{July 18, 2003. Revised February 11, 2004.}
\begin{document}
\begin{abstract}
We consider convex sets and functions over idempotent
semifields, like the max-plus semifield.
We show that if $\sK$ is a conditionally
complete idempotent semifield, with completion
$\kbar$, a convex function $\sK^n\to\kbar$ 
which is lower semi-continuous
in the order topology is the upper hull
of supporting functions defined as residuated differences
of affine functions. This result is proved using a separation
theorem for closed 
convex subsets of $\sK^n$, which extends earlier
results of Zimmermann, Samborski, and Shpiz. 
\end{abstract}
\maketitle
\section{Introduction}\label{sec-intro}
\begin{figure}[hbtp]
\begin{center}
\includegraphics[scale=0.55]{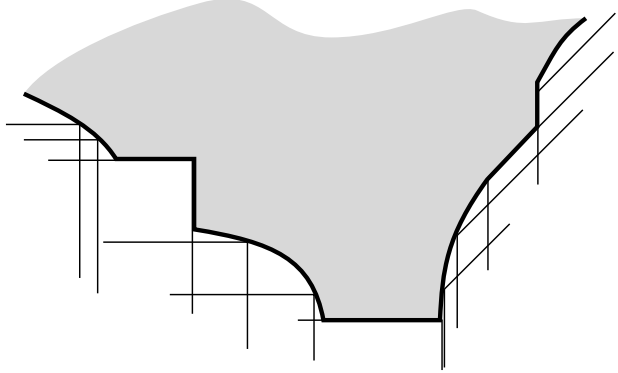}
\caption{A convex function over \(\rmax\) and its supporting half spaces}
\label{fig:fconvex}
\end{center}
\end{figure}

In this paper, we consider convex subsets
of semimodules over semirings with an idempotent addition,
like the max-plus semifield $\rmax$, which is the set $\RR\cup\{-\infty\}$,
with $(a,b)\mapsto \max(a,b)$ as addition, and $(a,b)\mapsto a+b$
as multiplication.  Convex subsets $C\subset \rmax^n$,
or \new{max-plus convex sets}, satisfy
\[
(x,y\in C,\; \alpha,\beta\in\rmax,\; \max(\alpha,\beta)=0)
\implies \max(\alpha+ x,\beta+y)\in C \enspace,
\]
where the operation ``$\max$''should be understood
componentwise, and where $\alpha +x=(\alpha
+x_{1},\ldots,\alpha +x_{n})$ for $x=(x_{1},\ldots,x_{n})$.
We say that a function $f:\rmax^n\to\rbar:=\RR\cup\{\pm\infty\}$
is \new{max-plus convex} if its epigraph is max-plus convex. 
An example of max-plus convex function is depicted in Figure~\ref{fig:fconvex}
(further explanations will be given in~\S\ref{sec-fconvex}).

Motivations to study semimodules and convex sets
over idempotent semirings arise from several
fields. First, semimodules over idempotent semirings,
which include as special cases sup-semilattices with a bottom element
(which are semimodules over the Boolean semiring), are natural
objects in lattice theory.
A second motivation arises from dynamic programming
and discrete optimization. Early results in this direction 
are due to Cuninghame-Green (see~\cite{cuning}),
Vorobyev~\cite{V67,V70}, Romanovski~\cite{R67a}, 
K.~Zimmermann~\cite{Zimmermann.K}.
The role of max-plus algebra in Hamilton-Jacobi
equations and quasi-classical asymptotics,
discovered by Maslov~\cite[Ch.~VII]{maslov73}
led to the development of an ``idempotent analysis'', 
by Kolokoltsov, Litvinov, Maslov, Samborski, Shpiz, 
and others (see~\cite{maslov92,maslovkolokoltsov95,litvinov00}
and the references therein).
A third motivation arises from
the algebraic approach of discrete event systems~\cite{bcoq}:
control problems for discrete event systems
are naturally expressed in terms
of invariant spaces~\cite{ccggq99}.

\begin{figure}[hbtp]
\begin{center}
\includegraphics[scale=0.65]{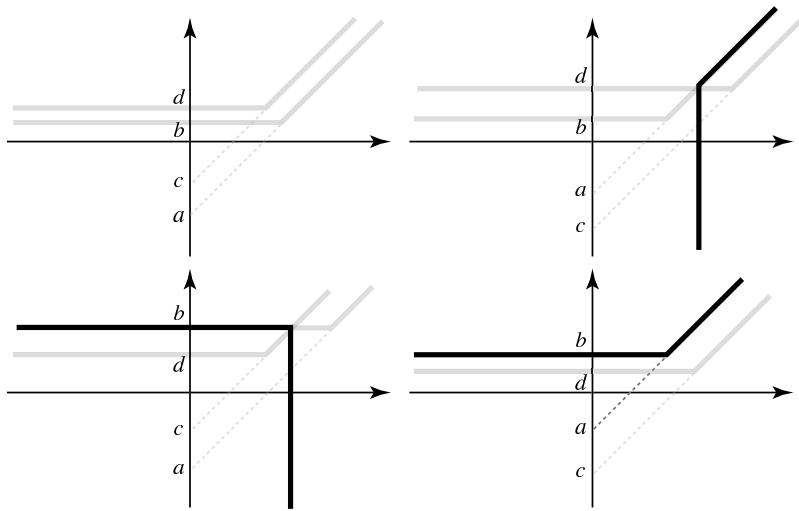}
\caption{The four generic differences of affine functions plots}
\label{fig:affine}
\end{center}
\end{figure}

Another motivation, directly related to the present work,
comes from abstract convex analysis~\cite{singer84,ACA,rubinov}:
a basic result of convex analysis
states that convex lower semi-continuous functions
are upper hulls of affine maps, which means precisely
that the set of convex functions is the max-plus 
(complete) semimodule generated by linear maps. In the theory of generalized
conjugacies, linear maps are replaced by a general family
of maps, and the set of convex functions is replaced
by a general semimodule. 
More precisely, given an abstract class of convex sets and functions,
a basic issue is to find a class of elementary functions 
with which convex sets and functions can be represented. This can be
formalized in terms of \new{$U$-convexity}~\cite{dk,ACA}.
If $X$ is a set and $U\subset \rbar^X$, a function
$f:X\rightarrow \rbar$ is called $U$-\emph{convex} if there exists a
subset $U'$ of $U$ such that 
\begin{equation}
f(x)=\sup_{u\in U'}u(x)\enspace,\qquad \forall x\in X\enspace.  \label{uconv}
\end{equation}
A subset $C$ of $X$ is said to be $U$-\emph{convex}~\cite{fan,ACA}
if for each $y\in X\setminus C$ we can find a map $u\in U$
such that 
\begin{align}
u(y) > \sup_{x\in C} u(x) \enspace .\label{uconv2}
\end{align}
In this paper, we address the problem of finding the set $U$
adapted to max-plus convex sets and functions.
The analogy with classical algebra suggests
to introduce \new{max-plus linear} functions:
\begin{align*}
\<a,x\>=\max_{1\leq i\leq n} (a_i+x_i) \enspace ,
\end{align*}
with $a=(a_i)\in \rmax^n$,
and \new{max-plus affine} functions,
which are of the form
\begin{align}
u(x)= \max(\<a,x\>,b)\enspace ,
\label{e-def-mplinear2}
\end{align}
where $b\in\rmax$. 
In the max-plus
case, we cannot take for $U$ the set of
affine or linear functions, 
because any sup of max-plus affine (resp.\ linear)
functions remains max-plus affine (resp.\ linear).
This is illustrated in the last (bottom
right) picture in Figure~\ref{fig:affine}, which shows
the graph of a generic affine function in dimension $1$
(the graph is the black broken line,
see Table~\ref{tab:1} in \S\ref{sec-fconvex} for details).
It is geometrically obvious that we cannot obtain the convex function
of Figure~\ref{fig:fconvex} as the sup of affine functions.
(Linear functions, however, lead to an
interesting theory if we consider
max-plus \new{concave} functions instead of max-plus convex
functions,  see Rubinov and Singer~\cite{singer00}, 
and downward sets instead of max-plus convex sets, 
see Mart\'\i nez-Legaz, Rubinov, and 
Singer~\cite{singer}.)

We show here that for max-plus convex functions, 
and more generally for convex functions over
conditionally complete idempotent semifields,
the
appropriate $U$ consists of \new{residuated differences}
of affine functions, which are of the form
$u\ominus u'$, where $u,u'$ are affine functions,
and $\ominus$ denotes the residuated law of the semiring addition,
defined in~\eqref{e-def-ominus} below. 
Theorem~\ref{tconverse} 
shows that lower semi-continuous convex functions
are precisely upper hulls of residuated differences of
affine functions, and Corollary~\ref{cor-Uconvex}
shows that the corresponding $U$-convex sets 
are precisely the closed convex sets.

As an illustration, in the case of the max-plus semiring,
in dimension $1$, there are 4 kinds of residuated
differences of affine functions, as shown in Figure~\ref{fig:affine},
and Table~\ref{tab:1}: one of these types 
consists of affine functions (bottom right, already discussed), 
another of these types consists only of the identically
$-\infty$ function (top left, not visible), whereas
the top right and bottom left plots yield new shapes, 
which yield ``supporting half-spaces'' for 
the convex function of Figure~\ref{fig:fconvex}.

The main device in the proof
of these results is a separation theorem for closed
convex sets (Theorem~\ref{th-1}).
We consider a convex subset $C$
of $\sK^n$, where $\sK$ is an
idempotent semifield 
that is conditionally complete for its natural
order. Then, we show that if $C$ is stable
under taking sups of (bounded) directed subsets
and infs of (bounded) filtered subsets,
or equivalently, if $C$ is closed in Birkhoff's order topology,
and if $y\in \sK^n\setminus C$,
there exists an affine hyperplane
\[
H=\bset x\in \sK^n\mset u(x) = u'(x) \eset  \enspace,
\]
with $u,u'$ as in~\eqref{e-def-mplinear2},
containing $C$ and not $y$.
When $\sK=\rmax$, Birkhoff's order
topology coincides with the usual one, and
we get a separation theorem for convex
subsets of $\rmax^n$ which are closed
in the usual sense. 
A key discrepancy,
by comparison with usual convex sets,
is that a two sided equation $u(x)=u'(x)$ is needed.
Theorem~\ref{th-1} extends or refines
earlier results by Zimmermann~\cite{zimmerman77},
Samborski and Shpiz~\cite{shpiz},
and by the three first authors~\cite{cgq02}.
Some metric assumptions on the semifield, which
were used in~\cite{zimmerman77}, are eliminated,
and the proof of Theorem~\ref{th-1} is in our view 
simpler (with a direct
geometric interpretation in terms of projections).
The method of~\cite{shpiz} only applies to the case
where the vector $y$ does not have entries equal
to the bottom element.
This restriction is removed in Theorem~\ref{th-1}
(see Example~\ref{ex-ncex} below for details). 
By comparison with~\cite{cgq02}, the difference
is that we work here in conditionally complete semifields
(without a top element), whereas the result of~\cite{cgq02}
applies to the case of complete semirings (which necessarily
have a top element). When the top element
is a coefficient of an affine equation defining an hyperplane,
the hyperplane need not be closed in the order topology,
and a key part of the proof of Theorem~\ref{th-1} is precisely
to eliminate the top element from the equation defining
the hyperplanes. In many applications, Birkhoff's
order topology is the natural one, so that the present
Theorem~\ref{th-1} gives a useful
refinement of the universal separation result of~\cite{cgq02}.

We finally point out additional references in which 
semimodules over idempotent semirings or related structures appear:
\cite{K65,Zimmermann.U,cao84,wagneur91,golan92,CGQ96a,CGQ97,nuclear,gondran02}. 
\section{Preliminaries}
\subsection{Ordered sets, residuation, idempotent semirings and semimodules}
In this section, we recall some basic notions
about partially ordered sets, residuation, idempotent semirings
and semimodules. See~\cite{birkhoff40,dubreil53,blyth72,cgq02}
for more details.  By \new{ordered set},
we will mean throughout the paper a set equipped
with a \new{partial} order. We say that an ordered set $(S,\leq)$
is \new{complete} if any subset $X\subset S$
has a least upper bound (denoted by $\supp X$).
In particular, $S$ has both a minimal (bottom) element 
$\minn S=\supp\emptyset$, and a maximal (top)
element $\maxx S=\supp S$. 
Since the greatest lower bound of a subset $X\subset S$ can
be defined by $\inff X=\supp\bset y\in S\mset y\leq x,\; \forall x\in X\eset$,
$S$ is a complete lattice. We shall also consider
the case where $S$ is only \new{conditionally complete},
which means that any subset of $S$
bounded from above has a least upper
bound and that any subset of $S$ bounded
from below has a greatest lower bound. 

If $(S,\leq)$ and $(T,\leq)$ are ordered
sets, we say that a map $f: S\to T$ is \new{residuated}
if there exists a map $f\sh: T\to S$ such that 
\begin{equation}
f(s) \leq t \iff s\leq f\sh(t) \enspace ,
\label{e-def-res}
\end{equation}
which means that for all $t\in T$, 
the set $\bset{s\in S}\mset{f(s)\leq t}\eset$
has a maximal element, $f\sh(t)$.
If $(X,\leq)$ is an ordered set,
we denote by $(X\op,\leqop)$
the \new{opposite} ordered
set, for which $x\leqop y\iff x\geq y$.
Due to the symmetry of the defining property~\eqref{e-def-res},
it is clear that if $f:S\to T$ is residuated,
then $f\sh: T\op\to S\op$ is also residuated.
When $S,T$ are complete ordered
sets, there is a simple
characterization of residuated maps.
We say that a map $f:S\to T$
\new{preserves arbitrary sups} if for all $U\subset S$,
$f(\supp U) = \supp f(U)$, where
$f(U)=\bset{f(x)}\mset{x\in U}\eset$.
In particular, when $U=\emptyset$, 
we get $f(\minn S)=\minn T$. 
One easily checks 
that
if $(S,\leq)$ and $(T,\leq)$ are complete ordered sets,
then, a map $f:S\to T$ is residuated if, and only if,
it preserves arbitrary sups
(see~\cite[Th.~5.2]{blyth72}, or~\cite[Th.~4.50]{bcoq}).
In particular, a residuated map $f$ is \new{isotone},
$x\leq y\implies f(x)\leq f(y)$, which, together
with~\eqref{e-def-res}, yields $f\circ f\sh\leq I_T$
and $f\sh\circ f\geq I_S$, where $I_T$ (resp.\ $I_S$)
denotes the identity map of $T$ (resp.\ $S$). This also implies that:
\begin{align}
f\circ f\sh\circ f=f \label{ffshf=f}, \qquad f\sh\circ f\circ f\sh=f\sh
\enspace .
\end{align}

We now apply these notions to idempotent semirings
and semimodules. 
Recall that a \new{semiring} is a set $\sS$ equipped
with an addition $\oplus$ and a multiplication $\otimes$, such that
$\sS$ is a commutative monoid for addition,
$\sS$ is a monoid for multiplication, multiplication
left and right distributes over addition, 
and the zero element of addition, $\zero$, is absorbing
for multiplication. We denote by $\unit$ the neutral
element of multiplication (unit).
We say that $\sS$ is \new{idempotent} when $a\oplus a=a$.
All the semirings considered in the sequel will be idempotent.
We shall adopt the usual conventions, and write for instance $ab$ instead
of $a\otimes b$.
An idempotent monoid $(S,\oplus,\zero)$
can be equipped with the \new{natural} order relation,
$a\leq b\Leftrightarrow a\oplus b=b$, for which $a\oplus b=a\vee b$,
and $\zero=\minn \sS$. We say that the semiring $\sS$ is
\new{complete} (resp.\ \new{conditionally complete})
if it is complete (resp.\ conditionally complete)
as a naturally ordered set,
and if for all $a\in \sS$, the left and right multiplications
operators, $\sS\to \sS$, $x\mapsto ax$, and $x\mapsto xa$,
respectively, preserve arbitrary sups
(resp.\ preserves sups of bounded from above sets).
An \new{idempotent semifield} is an idempotent semiring
whose nonzero elements are invertible. An idempotent
semifield $\sS$ cannot be complete, unless $\sS$
is the two-element \new{Boolean semifield},
$\{\zero,\unit\}$.
However, a conditionally complete semifield $\sS$
can be embedded in a complete semiring $\sbar$,
which is obtained by adjoining to $\sS$ a top
element, $\tau$, and setting
$a\oplus \tau =\tau$, $\zero\tau=\tau\zero=\zero$,
and $a\tau=\tau a=\tau$ for $a\neq\zero$.
Then, we say that $\sbar$ is the \new{completed
semiring} of $\sS$
($\sbar$ was called the \new{top-completion} of $\sS$ in~\cite{CGQ97},
and the \new{minimal completion} of $\sS$ in~\cite{akiansinger}).
For instance, the \new{max-plus semifield}
$\rmax$, defined in the introduction, can be embedded in the 
\new{completed max-plus semiring}
$\rmaxb$, whose set of elements is $\rbar$.

A (right) $\sS$-\new{semimodule} $X$
is a commutative monoid $(X,\oplus,\zero)$,
equipped with a map $X\times \sS\to X$, $(x,\lambda)
\to x\lambda$ (right action), that satisfies
\(x(\lambda \mu)=(x\lambda)\mu\),
\((x\oplus y)\lambda = x\lambda\oplus
y\lambda\), \(x(\lambda \oplus \mu) = x\lambda\oplus
x\mu\), \(x\zero = \zero\),
and \(x\unit=x\),
for all $x,y\in X$, $\lambda,\mu\in \sS$,
see~\cite{cgq02} for more details.
Since $(\sS,\oplus)$ is idempotent, 
$(X,\oplus)$ is idempotent, so that $\oplus$
coincides with the $\vee$ law for the natural
order of $X$. 
All the semimodules that we shall consider
will be right semimodules over idempotent semirings.
If $\sS$ is a complete semiring,
we shall say that a $\sS$-semimodule $X$
is \new{complete}
if it is complete as a naturally ordered set,
and if, for all $x\in X$ and
$\lambda\in \sS$,
the left and right multiplications,
$X\to X$, $x\mapsto x\lambda$,
and $\sS\to X$, $\mu\mapsto v\mu$,
respectively, preserve arbitrary sups.
We shall say that $V\subset X$ is
a \new{complete subsemimodule} of $X$
if $V$ is a subsemimodule of $X$ stable under arbitrary
sups. A basic example of semimodule 
over an idempotent semiring $\sS$ is
the \new{free semimodule} $\sS^n$,
or more generally the semimodule $\sS^I$ of
functions from an arbitrary set $I$ to $\sS$,
which is complete when $\sS$ is complete.
For $x\in \sS^I$ and $i\in I$, we denote, as
usual, by $x_i$ the $i$-th entry of $x$.

In a complete semimodule $X$, we define,
for all $x,y\in X$,
\begin{align*}
  x\lres y & = \maxx\bset{\lambda\in \sS}\mset{x\lambda \leq y}\eset \enspace,
\end{align*}
where we write $\maxx$ for the least upper bound to emphasize
the fact that the set has a top element.
In other words, $y\mapsto x\lres y, \; \sS\to \sS$ is
the residuated map of $\lambda \mapsto x\lambda, \; X\to X$. 
Specializing~\eqref{e-def-res}, we get
\begin{align}
\label{e-1}
x\lambda \leq y\iff \lambda \leq x\lres y 
\enspace .
\end{align}
For instance, when $\sS=\rmaxb$,
$(-\infty)\lres(-\infty) = 
(+\infty)\lres(+\infty)=+\infty$,
and $\mu \lres \nu = \nu - \mu$
if $(\mu,\nu)$ takes other values
($\sS$ being thought of as a semimodule over itself).
More generally, if $\sS$ is any complete semiring,
the law ``$\lres$'' of the semimodule $\sS^n$ can be
computed from the law ``$\lres$'' of $\sS$ by
\begin{align}
x\lres y= \inff_{1\leq i\leq n} x_i\lres y_i \enspace .
\label{e-def-lres}
\end{align}
Here, $\lres$ has a higher priority than $\wedge$,
so that the right hand side of~\eqref{e-def-lres}
reads $\inff_{1\leq i\leq n} (x_i\lres y_i)$.
If the addition of $\sS$ distributes over arbitrary infs
(this is the case in particular if $\sS$ is a semifield,
or a completed semifield, see~\cite[Ch.~12, Th.~25]{birkhoff40}),
for all $\lambda\in \sS$, the translation
by $\lambda$, $\mu\mapsto \lambda \oplus \mu$, 
defines a residuated map $\sS\op\to\sS\op$,
and we set: 
\begin{align}
\nu\ominus \lambda = \minn\bset \mu \mset \lambda \oplus\mu \geq \nu \eset 
\enspace , \label{e-def-ominus}
\end{align}
where we write $\minn$ for the greatest lower bound to emphasize
the fact that the set has a bottom element.
When $\sS=\rmaxb$, we have (see e.g.\ \cite{bcoq,singer91}):
\[
\nu\ominus \mu= \begin{cases}
\nu & \mrm{ if}\;\nu>\mu,\\
-\infty & \mrm{otherwise.}
\end{cases}
\]
Dualizing the definition~\eqref{e-def-res} of residuated maps, 
we get:
\begin{align}
\label{e-1dual}
\lambda\oplus\mu \geq \nu \iff \lambda \geq \nu\ominus \mu 
\enspace .
\end{align}

\subsection{Separation theorem for complete convex sets}
We next recall the general separation
theorem of~\cite{praha,cgq02}.
By \new{complete semimodule}, we mean throughout
the section a complete semimodule
over a complete idempotent semiring $\sS$.
Let $V$ denote a complete subsemimodule of a complete semimodule $X$.
We call \new{canonical projector} onto $V$ the map
\[
P_V: X\to V,\quad P_V(x) = \maxx\bset{v\in V}\mset{v\leq x}\eset	\]
(the least upper bound of $\bset{v\in V}\mset{v\leq x}\eset$ belongs
to the set because $V$ is complete).
Thus, $P_V$ is the residuated map of the canonical
injection $i_V: V\to X$, $P_V$ is surjective,
and $P_V=P_V^2$.
If $\{w_\ell\}_{\ell\in L}\subset X$ is an arbitrary
family, we set 
\[
\bigoplus_{\ell\in L} w_\ell :=\supp \bset w_\ell \mset \ell\in \ELL\eset
\enspace .
\]
We say that $W$ is a \new{generating family}
of a complete subsemimodule $V$ if any
element $v\in V$ can be written as 
$v=\bigoplus_{w\in W} w \lambda_w$,
for some $\lambda_w\in \sS$.
If $V$ is a complete subsemimodule of $X$ with generating
family $W$, then
\begin{align}
P_V(x) = \bigoplus_{w\in W} w (w\lres x) \enspace ,
\label{e-proj}
\end{align}
see~\cite[Th.~5]{cgq02}.
\begin{theorem}[Universal Separation Theorem,~{\cite[Th.~8]{cgq02}}]
\label{th-separ}
Let $V\subset X$ denote a complete subsemimodule, 
and let $y\in X\setminus V$. Then, the set
\begin{align}
H= \bset{x\in X}\mset{x\lres P_V(y)=x\lres y}\eset 
\label{e-def-seph}
\end{align}
contains $V$ and not $y$.
\end{theorem}
Seeing $x\lres y$ as a ``scalar product'',
$H$ can be seen as the ``hyperplane''
of vectors $x$ ``orthogonal'' to
$(y,P_V(y))$. As shown in~\cite{cgq02}, the ``hyperplane''
$H$ is a complete subsemimodule
of $X$, even if it is defined by a nonlinear equation.
In order to give a linear defining equation
for this hyperplane, we have to make
additional assumptions on the semiring $\sS$.
In this paper,  we shall assume that $\sS=\kbar$
is the completed semiring of 
a conditionally complete idempotent semifield $\sK$.
Consider the semimodule of functions $X=\kbar^I$.
When $x=(x_i)_{i\in I}, y=(y_i)_{i \in I} \in \kbar^I$, 
we define 
\begin{align}
\<y,x\>=\bigoplus_{i\in I} y_i x_i 
\label{e-mpdotprod}
\end{align}
and 
\begin{align*}
\ldual{x}= \maxx\bset{y\in \kbar^I}\mset{\<y,x\>\leq \unit}\eset
\end{align*}
that is, 
\begin{align}
(\ldual{x})_i =x_i \lres \unit\enspace .
\label{three13}
\end{align}
(For instance, when $\sK=\rmax$, 
$(\ldual{x})_i = - x_i$.) We have
$\ldual{(x\lres y)}= \<\ldual{y}, x\>$,
and $\lambda \mapsto \ldual{\lambda}$ is bijective
$\kbar\to\kbar$, which allows us to write $H$ linearly:
\begin{align}
\label{e-linear}
H= \bset{x\in \kbar^I}\mset{\< \ldual{P_V(y)}, x\> = \<\ldual{y}, x\>}\eset 
\end{align}
(see~\cite{cgq02} for generalizations to more general semirings,
called \new{reflexive} semirings).

Theorem~\ref{th-separ} yields a separation
result for convex sets as a corollary. 
We recall that a subset $C$ of a complete semimodule $X$
over a complete semiring $\sS$ is \new{convex}~\cite{zimmerman77,zimmermann}
(resp.\ \new{complete convex}~\cite{cgq02})
if for all finite (resp.\ arbitrary)
families $\{x_\ell\}_{\ell\in \ELL}\subset C$ and $\{\alpha_\ell\}_{\ell\in \ELL}
\subset \sK$, such that $\bigoplus_{\ell\in \ELL} \alpha_\ell =\unit$,
we have that $\bigoplus_{\ell\in \ELL} x_\ell\alpha_\ell \in C$. 
For example, every subsemimodule of $X$ is convex, 
and every complete subsemimodule of $X$ is complete convex.
\begin{corollary}[Separation Theorem for Complete Convex Sets, {\cite[Cor.~15]{cgq02}}]\label{cor-convex}
If $C$ is a complete convex subset of a complete semimodule $X$,
and if $y\in X\setminus C$, then the set 
\begin{align}
 H= \bset{x\in X}\mset{x\lres y \wedge \unit=x\lres Q_C(y) \wedge \nu_C(y)}\eset\label{e-formul0}
\end{align}
with
\begin{align}
\nu_C(y) &=\supp_{v\in C}(v\lres y\wedge \unit)\qquad \text{and}\qquad Q_C(y) = \supp_{v\in C}v(v\lres y\wedge \unit)   \;,\label{e-formul2}
\end{align}
contains $C$ and not $y$. 
\end{corollary}
Recall our convention explained after Equation~\ref{e-def-lres},
that $\lres$ has a higher priority than $\wedge$, so that for instance 
$v\lres y \wedge \unit=(v\lres y)\wedge \unit$.
When $\sS=\kbar$ is a completed idempotent semifield,
and $X=\kbar^I$, $H$ can be rewritten linearly:
\begin{align}
\label{e-separh-convex}
 H= \bset{x\in \kbar^I}\mset{\< \ldual{y}, x\>\oplus\unit = \<\ldual{Q_C(y)}, x\>\oplus \ldual{\nu_C(y)}}\eset  \enspace .
\end{align}
\begin{remark}\label{rk-ineq}
Since $Q_C(y)\leq y$, and $\nu_C(y)\leq\unit$, 
we have $\ldual{y}\leq \ldual{Q_C(y)}$
and $\unit=\ldual{\unit}\leq \ldual{\nu_C(y)}$,
and hence, by definition of the natural order $\leq$,
we can write equivalently
$H$ as
\[
H= 
\bset{x\in \kbar^I}\mset{\< \ldual{y}, x\>\oplus\unit 
\geq  \<\ldual{Q_C(y)}, x\>\oplus \ldual{\nu_C(y)}}\eset  \enspace .
\]
The same remark applies, mutatis mutandis, to~\eqref{e-def-seph},
~\eqref{e-linear}, and ~\eqref{e-formul0}.
\end{remark}
\subsection{Geometric interpretation}
We complement the results of~\cite{cgq02} 
by giving a geometric interpretation 
to the vector $Q_C(y)$ and scalar $\nu_C(y)$
which define the separating hyperplane $H$.
If $C$ is any subset of $X$, we call \new{shadow of \(C\)},
denoted by $\Down(C)$, the set of linear combinations
\[
\bigoplus_{\ell\in \ELL}x_\ell\lambda_\ell,\;\text{ with }\; x_\ell\in C, \;\lambda_\ell\in \sS,\;
\lambda_\ell\leq\unit\enspace,
\text{ and } \ELL \text{ a possibly infinite set.}
\]
We also denote by
\[
\Up(C)=\bset{z\in C}\mset{\exists v\in C,\; z\geq v}\eset
\]
the \new{upper set} generated by $C$. 
The term ``shadow'' can be interpreted geometrically: when 
for instance $C\subset \rmaxb^2$,
$\Down(C)$ is the shadow of $C$ if the sun light comes
from the top-right corner of the plane, see Figure~\ref{fig:shadowup}
and Example~\ref{ex-sh} below.
\begin{theorem}[Projection onto $\Down(C)$ and $C$]\label{theo-qc}
If $C$ is a complete convex subset of a complete
semimodule $X$, then, for all $y \in X$,
\begin{align}
Q_C(y)=\maxx \bset{z\in \Down(C)}\mset{z\leq y}\eset 
\; .
\label{e-downbarqc}
\end{align}
If $y\in \Up(C)$,
\begin{align}
Q_C(y)=\maxx \bset{z\in C}\mset{z\leq y}\eset ,\qquad \text{and}\qquad
\nu_C(y) =\unit 
\; .
\label{e-newqc}
\end{align}
If \(\nu_C(y)\) is invertible, \(Q_C(y)(\nu_C(y))^{-1}\) belongs to $C$.
\end{theorem}
Thus, Theorem~\ref{theo-qc} shows that $Q_C$ is a projector
which sends $X$ to $\Down(C)$, and $\Up(C)$ to $C$. Moreover, when  \(\nu_C(y)\) is invertible, \(Q_C(y)(\nu_C(y))^{-1}\) can be considered
as the projection of $y$ onto $C$.
\begin{proof}
Since $v\lres y\wedge\unit \leq \unit$,
\begin{align}
Q_C(y)=\bigoplus_{v\in C} v (v\lres y\wedge\unit)
\in \Down(C)
\enspace .\label{e-XtoC}
\end{align}
If $y\in \Up(C)$, we have $v\leq y$ for some $v\in C$,
hence, $v\lres y\geq \unit$ (by~\eqref{e-1}), which implies that
$\nu_C(y) \geq v\lres y\wedge \unit =\unit$.
Since $\nu_C(y)\leq \unit$ holds trivially,
we have proved that $\nu_C(y)=\unit$, so that
\begin{align}
y\in \Up(C)\implies
Q_C(y)\in C
\text{ and } \nu_C(y)= \unit \enspace .
\label{e-UptoC}
\end{align}
Consider now any element $z\in \Down(C)$,
$z=\bigoplus_{\ell\in \ELL} v_\ell\lambda_\ell$,
with $v_\ell\in C$, $\lambda_\ell\in \sS$,
$\lambda_\ell\leq \unit$, and assume that $z\leq y$.
Then, $v_\ell\lambda_\ell\leq y$,
so that $\lambda_\ell \leq v_\ell\lres y$
(by~\eqref{e-1}), and since $\lambda_\ell\leq \unit$,
$Q_C(y) \geq v_\ell (v_\ell\lres  y \wedge \unit)
\geq v_\ell (\lambda_\ell\wedge \unit)
=v_\ell \lambda_\ell$. Summing
over all $\ell\in L$, we get $Q_C(y)\geq z$.
Together with~\eqref{e-XtoC}, this shows~\eqref{e-downbarqc}.
Since we also proved~\eqref{e-UptoC}, this shows a fortiori~\eqref{e-newqc}.

Finally, if $\nu_C(y)$ is invertible, 
we see from~\eqref{e-XtoC}
that \(Q_C(y)(\nu_C(y))^{-1}\) is of the form $\bigoplus_{v\in C}v\lambda_{v}$ with $\bigoplus_{v\in C}\lambda_{v}=\unit$, hence \(Q_C(y)(\nu_C(y))^{-1}\) belongs to $C$.
\end{proof}
\begin{example}\label{ex-sh}
In Figure~\ref{fig:shadowup}, the convex~\(C\) generated by three points \((a,b,c)\) in \(\rmaxb^{2}\) is displayed, together with its shadow and upper set. The cases of \(y\) belonging to \(\Up(C)\) and of \(y\notin \Up(C)\) are illustrated.
\end{example}
\begin{figure}
\begin{center}
\includegraphics{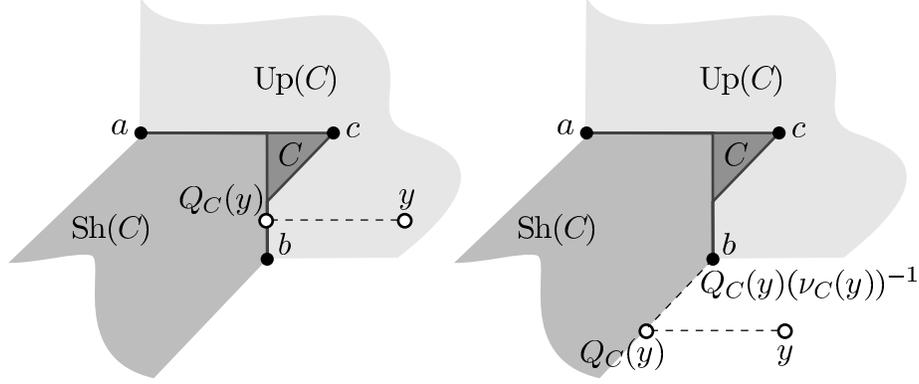}
\caption{Projections}
\label{fig:shadowup}
\end{center}
\end{figure}
\begin{remark}
When $\sS=\kbar$ is a completed idempotent semifield, and
$C$ is complete and convex, then
\begin{align}
\Down(C)=\bset x \lambda \mset \lambda\in \sK,\;\lambda\leq\unit,\;x\in C\eset\enspace.
\label{e-down-semifield}
\end{align}
Indeed, let $\Down'(C)$ denote the set in the right hand side of~\eqref{e-down-semifield}. The inclusion $\Down'(C)\subset \Down(C)$
is trivial. To show the other
inclusion, take any $z\in\Down(C)$,
which can be written as a linear combination $z=\bigoplus_{\ell\in \ELL}x_\ell\lambda_\ell$,
for some $x_\ell\in C$, $\lambda_\ell\in \sS$, $\lambda_\ell\leq \unit$, 
with $\ELL$ a possibly infinite set.
When $z=\zero$, $z\in \Down'(C)$ trivially.
When $z\neq \zero$, $\lambda_\ell\neq\zero$
for some $\ell$, so that $\mu:=\bigoplus_{\ell\in \ELL}\lambda_\ell\neq\zero$,
and since $\mu\leq\unit$ and $\sS$ is a completed idempotent
semifield, $\mu$ is invertible.
Writing $z=y \mu$, and observing that
$y=\bigoplus_{\ell\in \ELL}x_\ell\lambda_\ell\mu^{-1}$
belongs to $C$ because $C$ is complete and convex,
we see that $z\in \Down'(C)$.\qed
\end{remark}
\begin{example}\label{cex}
To illustrate the previous results,
consider the convex set $C\subset \rmaxb^2$ generated
by the two points $(0,-\infty)$ and $(2,3)$.
Thus, $C$ is the set of points of the form
$(\max(\alpha,\beta+2), \beta+3)$,
with $\max(\alpha,\beta)=0$. Since $C$ is generated
by a finite number of points of $\rmaxb^2$,
$C$ is complete convex. The set $C$ is the broken
dark segment between the points $(0,-\infty)$
and $(2,3)$, in Figure~\ref{fig-todo}.
In order to represent points with $-\infty$
coordinates, we use exponential coordinates
in Figure~\ref{fig-todo}, that is, the point
$(z_1,z_2)\in \rmax^2$ is
represented by the point of the positive
quadrant of coordinates $(\exp(z_1),\exp(z_2))$.
Consider now $y=(1,-k)$,
for any $k\geq 0$, and let us separate $y$ from
$C$ using Corollary~\ref{cor-convex}. Since $(0,-\infty)\leq
y$, $y\in \Up(C)$, and we get from~\eqref{e-newqc} that $\nu_{C}(y)=0$.
One also easily checks that $Q_C(y)=(0,-k)$. 
When $k\neq+\infty$, the separating hyperplane $H$ of~\eqref{e-separh-convex}
becomes:
\begin{align}
\label{e-sep-h}
H= \bset{ x\in \rmaxbv 2}\mset{\max(-1+ x_1,k+x_2,0)=\max(x_1, k+x_2,0)}\eset
\enspace .
\end{align}
The point $y=(1,0)$, together
with $Q_C(y)$ and the separating hyperplane $H$
(light grey zone) are depicted at the left of Figure~\ref{fig-todo}. 
When $k=\infty$, is it easily checked
that the separating hyperplane is
the union of the half space 
$\rbar\times (\RR\cup\{\infty\})$,
and of the interval $[-\infty, 0]\times \{-\infty\}$.
Unlike in the case of a finite $k$,
$H$ is not closed for the usual topology,
which implies that the max-plus linear forms which
define $H$ are {\em not} continuous for the usual topology.
\begin{figure}[hbtp]
\begin{center}
\includegraphics{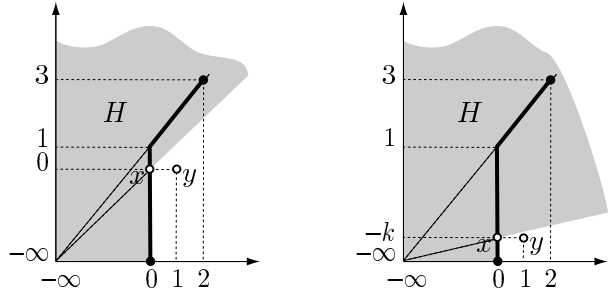}
\caption{Separating a point from a convex set}
\label{fig-todo}
\end{center}
\end{figure}
\end{example}
\subsection{Closed convex sets in the order topology}
We next recall some basic facts about Birkhoff's order topology~\cite[Ch.~10, \S~9]{birkhoff40}, and establish some properties of closed convex sets.
See~\cite{gierzETAL,akiansinger} for more background
on topologies on lattices and lattices ordered groups.

Recall that a nonempty ordered set $D$ is \new{directed}
if any finite subset of $D$ has an upper bound in $D$,
and that a nonempty ordered set $F$ is \new{filtered}
if any finite subset of $F$ has a lower bound in $F$.
\begin{definition}
We say that a subset $X$ of a conditionally complete
ordered set $S$ is \new{stable under directed sups} (resp.\ \new{stable under filtered infs})
if for all directed (resp filtered) subsets $D\subset X$ (resp.\ $F\subset X$)
bounded from above (resp.\ below), $\supp D\in X$
(resp.\ $\inff F\in X$).
\end{definition}
When $S=\sK^n$, 
where $\mathcal{K}$ is a conditionally complete idempotent semifield,
the condition that $F$ is bounded from below can be dispensed with,
since any $F$ is bounded from below by $\zero$.
Recall that a \new{net}
with values in 
a conditionally complete ordered set $S$
is a family $(x_\ell)_{\ell\in \ELL}\subset S$
indexed by elements of a directed set $(\ELL,\leq)$.
We say that a net $(x_\ell)_{l\in\ELL}\in S$ bounded from above
and from below \new{order converges} to $x\in S$, if
$x=\limsup_{\ell\in \ELL} x_\ell=\liminf_{\ell\in \ELL} x_\ell$,
where $\limsup_{\ell\in \ELL} x_\ell:=\inf_{\ell\in \ELL}\sup_{m\geq \ell} x_m$,
and $\liminf_{\ell\in \ELL}x_\ell:=\sup_{\ell\in \ELL}\inf_{m\geq \ell} x_m$.
We say that $X\subset S$ is \new{order-closed}
if for all nets $(x_\ell)_{\ell\in \ELL}\subset X$ order
converging to some $x\in S$, $x\in X$. The set
$o(S)$ of order-closed subsets of $S$ defines the Birkhoff's 
\new{order topology}.
In particular, if $D$ (resp.\ $F$) is a directed (resp.\ filtered)
subset of $X$, $\{x\}_{x\in D}$ (resp.\ $\{x\}_{x\in F^{\mrm{op}}}$)
is a net which order converges
to $\supp D$ (resp.\ $\inff F$), so that any order
closed set is stable under directed sups and filtered
infs.
We warn the reader that a net which is order
convergent is convergent for the order
topology, but that the converse need not hold,
see~\cite[Ch.~10, \S~9]{birkhoff40}.
However, both notions coincide when $S=\sK^n$ if $\sK$ is a conditionally
complete semifield which is a continuous lattice~\cite{akiansinger}.
When $S=\rmaxv n$,
the order topology is the usual topology
on $(\RR\cup\{-\infty\})^n$.

The following result applies in particular to convex
subsets of semimodules.
\begin{proposition}
A subset $C\subset S$ stable under finite sups is closed
for the order topology if and only if it is
stable under directed sups and filtered infs.
\end{proposition}
\begin{proof}
Assume that $C$ is stable under directed sups and
filtered infs, and let 
$\{x_\ell\}_{\ell\in \ELL}\subset C$ denote a net
order converging to $x\in S$. We have
$x=\inff_{\ell\in \ELL} \bar x_\ell$, where
$\bar x_\ell=\bigoplus_{m\geq \ell} x_m$.
Let $D_\ell$ denote the set of finite subsets of
$\bset{m\in \ELL}\mset{m\geq \ell}\eset$, 
and for all $J\in D_\ell$, define $x_J=\bigoplus_{m\in J} x_m$.
Since $C$ is stable under finite sups, $x_J\in C$.
Since $\bset{x_J}\mset{J\in D_\ell}\eset$ is directed and bounded from above,
and since $C$ is stable under directed sups, 
$\bar x_\ell =\bigoplus_{J\in D_\ell}\bigoplus_{m\in J} x_m
=\bigoplus_{J\in D_\ell} x_J\in C$. Since $\bset{\bar x_\ell}\mset \ell\in \ELL\eset$
is filtered and bounded from below, and since $C$ is stable under
filtered infs, $x=\inff_{\ell\in \ELL}\bar x_\ell\in C$,
which shows that $C$ is closed for the order topology.
This shows the ``if'' part of the result.
Conversely, if $D$ (resp.\ $F$) is a directed (resp.\ filtered) subset of $X$,
then $\{x\}_{x\in D}$ (resp.\ $\{x\}_{x\in F\op}$) is a net which order converges to 
$\supp D$ (resp.\ $\inff F$), so that any order
closed set is stable under directed sups and filtered infs.
\end{proof}

We shall use repeatedly the following lemma in the sequel.
\begin{lemma}[See.~{\cite[Ch.~13, Th.~26]{birkhoff40}}]
\label{lem-birkhoff}
If $\sK$ is a conditionally complete semifield,
if $x_\ell\in \sK$ order converges to $x\in\sK$,
and $y_\ell\in \sK$ order converges
to $y\in\sK$, then
$x_\ell\wedge y_\ell$ order
converges to $x\wedge y$,
$x_\ell\oplus y_\ell$ order
converges to $x\oplus y$,
and $x_\ell y_\ell$ order
converges to $x y$.
\end{lemma}
In fact, the result of~\cite{birkhoff40} is stated
only for elements of $\sKp$,
but the extension to $\sK$ is plain,
since $x\zero =\zero x=x,\,x\oplus \zero =\zero \oplus x=x$ and $x\wedge
\zero =\zero \wedge x=\zero $ for all $x\in \sK$.
However, Lemma~\ref{lem-birkhoff}
does {\em not} extend to $\kbar$: 
for instance, in $\rmaxb$, $(-\ell)_{\ell\in \NN}$
order converges to $-\infty$, but $((+\infty)+(-\ell))_{\ell\in \NN}$,
which is the constant sequence with value $+\infty$,
does not order converge to $(+\infty)+(-\infty)=-\infty$.
This is precisely why the separating hyperplane
provided by the universal separation theorem
need not be closed, see Example~\ref{cex} above.

\begin{corollary}\label{cor-order}
If $v\in \sK^n$, $w\in\sK^n\setminus\{\zero\}$,
if $x_\ell\in \sK^n$ order converges to $x\in\sK^n$,
and if $\lambda_\ell\in\sK$ order converges to $\lambda\in\sK$,
then, $\<v,x_\ell\>$ order converges
to $\<v,x\>$, $w\lres x_\ell$ order converges
to $w\lres x$, and $v\lambda_\ell$ order converges
to $v\lambda$.
\end{corollary}
\begin{proof}
By~\eqref{e-mpdotprod}, 
$\<v,y\>= \bigoplus_{1\leq i\leq n} v_i y_i$,
and by~\eqref{e-def-lres},
$w\lres y=
\inff_{i\in I} w_i^{-1} y_i$, 
where $I=\bset 1\leq i\leq n\mset w_i\neq\zero\eset\neq\emptyset$,
so the corollary follows from Lemma~\ref{lem-birkhoff}.
\end{proof}
We shall need the following basic property:
\begin{lemma}\label{l-1}
If $C$ is a convex subset (resp.\ a subsemimodule) of $\sK^n$,
then, its closure for the order topology is a convex subset
(resp.\ a subsemimodule) of $\sK^n$.
\end{lemma}
\begin{proof}
We derive this from Lemma~\ref{lem-birkhoff}
(the only unusual point is that the order convergence need not 
coincide with the convergence for the order topology).
Assume that $C$ is convex
(the case when $C$ is a semimodule is similar).
Recall that if $f$ is a continuous self-map
of a topological space $X$, then 
$f(\clo(Y))\subset \clo f(Y)$
holds for all $Y\subset X$, where $\clo(\cdot)$ denotes the
closure of a subset of $X$. 
Fix $\alpha,\beta\in \sK$
such that $\alpha\oplus \beta=\unit$, and 
consider $\psi: \sK^n\times \sK^n\to\sK^n$,
$\psi(x,y)=x\alpha \oplus y\beta$. 
We claim that for all $x\in \sK^n$, the map
$\psi(x,\cdot)$ is continuous in the order
topology. Indeed, let $A$ denote a subset
of $\sK^n$ that is closed in the order topology,
and let us show that the pre-image by $\psi(x,\cdot)$
of $A$, $A'=\bset y\in\sK^n\mset x\alpha\oplus y\beta \in A\eset$,
is also closed in the order topology.
If $\{y_\ell\}_{\ell\in \ELL}$ is any net
in $A'$ converging to some $y\in\sK^n$, we have
$x\alpha\oplus y_\ell\beta \in A$,
for all $\ell\in \ELL$, and it follows from Lemma~\ref{lem-birkhoff}
that $x\alpha\oplus y_\ell\beta$ order converges
to $x\alpha\oplus y\beta$. Since $A$ is closed
in the order topology, $x\alpha\oplus y\beta\in A$,
so $y\in A'$, which shows that $A'$ is closed
in the order topology. Thus, $\psi(x,\cdot)$
is continuous, and so $\psi(x,\clo(C))\subset \clo(\psi(x,C))$.
Since $C$ is stable under convex combinations, $\psi(x,C)\subset C$,
hence, $\psi(x,\clo(C))\subset \clo(C)$.
Pick now any $y\in \clo(C)$. Since
$\psi(x,y)\in \clo(C)$, for all $x\in C$,
and since $\psi(\cdot,y)$ is continuous,
$\psi(\clo(C),y)\subset \clo(\psi(C,y))\subset \clo(C)$.
Since this holds for all $y\in \clo(C)$, we have
shown that 
$\psi(\clo(C),\clo(C))\subset \clo(C)$,
i.e., $\clo(C)$ is stable under convex combinations.
\end{proof}
We conclude this section with properties
which hold more generally in semimodules of functions.
For all $C\subset \sK^I$, we denote by $\bar C\subset\kbar^I$
the set of arbitrary convex combinations of elements
of $C$: 
\begin{align}
\bar C=\bset \bigoplus_{\ell\in \ELL} v_\ell \lambda_\ell \mset 
\{v_\ell\}_{\ell\in \ELL} \subset C,\; 
\{\lambda_\ell\}_{\ell\in \ELL} \subset\sK,\;
\bigoplus_{\ell\in \ELL} \lambda_\ell=\unit\eset
\label{e-def-barc}
\end{align}
($\ELL$ denotes an arbitrary - possibly infinite - index set).
\begin{proposition}\label{prop-cbar1}
If $C$ is a convex subset of $\sK^I$ which is 
closed in the order topology, then
\begin{align}
\bar C\cap \sK^I= C \enspace.\label{e-fond} 
\end{align}
\end{proposition}
\begin{proof}
Consider an element
$v=\bigoplus_{\ell\in \ELL} v_\ell \lambda_\ell \in \bar C\cap \sK^I$,
with $\bigoplus_{\ell\in \ELL}\lambda_\ell=\unit$.
Assume, without loss of generality, that $\lambda_\ell\neq\zero$,
for all $\ell\in \ELL$.
Let $D$ denote the set of finite subsets of $\ELL$, and
for all $J\in D$, let $v_J=\bigoplus_{\ell\in J} v_\ell\lambda_\ell$,
and $\lambda_J= \bigoplus_{\ell\in J} \lambda_\ell$. By construction,
the net $\{v_J\}_{J\in D}$ order converges to $v$,
and the net $\{\lambda_J\}_{J\in D}$ order converges
to $\unit$. Hence, by Lemma~\ref{lem-birkhoff}, 
$v_J\lambda_J^{-1}$ order converges to $v$. 
But $v_J\lambda_J^{-1}=\bigoplus_{\ell\in J}v_\ell\lambda_\ell\lambda_{J}^{-1}\in C$,
and since $C$ is closed for the order topology, 
$v\in C$, which shows~\eqref{e-fond}. 
\end{proof}
When $C$ is a semimodule, the condition that $C$ is stable
under filtered infs, which is implied by the condition
that $C$ is closed in the order topology,
can be dispensed with.
\begin{proposition}\label{prop-cbar2}
If $C$ is a subsemimodule of $\sK^I$ which is 
stable under directed sups, ~\eqref{e-fond} holds.
\end{proposition}
\begin{proof}
Any element $v\in \bar C$ can be
written as $v=\bigoplus_{\ell\in \ELL} v_\ell$, for
some $\{v_\ell\}_{\ell\in \ELL}\subset C$. 
Setting $v_J=\bigoplus_{\ell\in J} v_\ell\in C$,
we get $v=\bigoplus_{J\in D} v_J$,
and we only need to know that $C$ is stable under directed sups
to conclude that $v\in C$.
\end{proof}
The following example shows that we
cannot derive Proposition~\ref{prop-cbar2}
from Proposition~\ref{prop-cbar1}.
\begin{example}
The set $C=\{(-\infty,-\infty)\}\cup (\RR\times \RR)$ is a subsemimodule
of $\rmax^2$, which is stable under directed
sups, but not stable under filtered infs (for instance
$\inff\bset (0,-\ell)\mset \ell\in\NN\eset=(0,-\infty)\not\in C$),
and hence not closed in the order topology.
\end{example}
\section{Separation theorems for closed convex sets}
\label{newsepar}
We saw in Example~\ref{cex} that, when $\sK=\rmax$,
the separating set~\eqref{e-separh-convex} given by the universal
separation theorem need not be closed for the usual topology.
In this section, we refine the universal
separation theorem in order to separate
a point from a {\em closed}
convex set by a {\em closed} hyperplane.

From now on, we assume that $\sK$
is a conditionally complete idempotent semifield,
whose completed semiring is denoted by $\kbar$.
\subsection{Separation of closed convex subsets of $\protect\mathcal{K}^I$}
As a preparation for the main result of \S\ref{newsepar}
(Theorem~\ref{th-1} below), we derive from
Corollary~\ref{cor-convex} a separation result
for order closed convex sets $C$ and elements $y\in X\setminus C$
of the semimodule of functions $X=\sK^I$,
satisfying an archimedean condition. This archimedean
condition will be suppressed in Theorem~\ref{th-1},
assuming that $I$ is finite.

\begin{definition}
\label{dhyperpl}
We call \new{affine hyperplane} of $\sK^I$ a subset
of $\sK^I$ of the form
\begin{align}
\label{hyperplane}
H= \bset{v\in \sK^I}\mset{\< w', x\>\oplus d' = \< w'', x\>\oplus d''}\eset  \enspace ,
\end{align}
with $w',w''\in \sK^I$, and $d',d'\in \sK$.
We shall say
that $H$ is a \new{linear hyperplane} if $d'=d''=\zero$.
\end{definition}
(When $I$ is infinite, $\<w',x\>$ and $\<w'',x\>$ may be
equal to $\top$.)
\begin{remark}
We have already encountered ``hyperplanes'' of $\sK^{I}$
of the above form.
Indeed, $H\cap \sK^{I}$, with $H$ of~\eqref{e-separh-convex},
is of the form~\eqref{hyperplane}, with 
\begin{equation}
w'=\ldual{y},\;d'=\unit,\;w''=\ldual{Q_{C}(y)},\;d''=\ldual{\nu _{C}(y)}
\enspace .  \label{two26}
\end{equation}
The main point in Definition~\ref{dhyperpl} is 
the requirements that $w',w''\in \sK^{I}$, and $d',d''\in \sK$
which need not be satisfied in (\ref{two26});
indeed, for $y=(y_{i})\in \sK^{I}$ having a
coordinate $y_{i_{0}}=\zero$, by \eqref{three13}, we
have $\ldual{y}_i=\top\kbar $,
so that $\ldual{y}\not\in \sK^{I}$
(see e.g.\ Example~\ref{cex}).
\end{remark}
Given $y\in \sK^I\setminus C$,
the question is whether we can find an
affine hyperplane of $\sK^I$
containing $C$ and not $y$. 
We shall need the following \new{Archimedean type assumption}
on $C$ and $y$:
\begin{align*}
(A):\qquad 
\forall v\in C,
\exists \lambda\in \sK\setminus\{\zero\}, 
v\lambda \leq y \enspace.
\end{align*}
For all $y\in \sK^I$ and $C\subset\sK^I$, 
define 
\[
\support y = \bset i\in I\mset y_i\neq \zero\eset \enspace ,
\qquad 
\support C= \bigcup_{v\in C} \support v \enspace .
\]
One readily checks that
$y\lambda \leq y' $ for some $\lambda\in \sK\setminus\{\zero\}$
implies that $\support y\subset \support y'$, 
and that when $I$ is finite,
the converse implication holds
(indeed, if 
$\support y\subset \support\,y'$, take
any $\lambda\in\sK$
smaller than $y\lres y'=\wedge_{i\in I}y_{i}\lres y_{i}'$,
a quantity which is in $\kbar\setminus \{\varepsilon \}$ when $I$ is
finite).
Thus, Assumption~(A) implies that 
\begin{align}
\support y \supset \support C \enspace ,
\label{e-sup=}
\end{align}
and it is equivalent to~\eqref{e-sup=} when $I$ is finite.
\begin{proposition}\label{th-n2}
Let $C$ be a convex subset of $\sK^I$,
and $y\in \sK^I\setminus C$.
Assume that $C$ is closed for the order topology of $\sK^I$,
and that Assumption (A) is satisfied.
Then, there is an affine
hyperplane of $\sK^I$ which contains $C$ and not $y$. 
\end{proposition}
\begin{remark}\label{rk-wpgeqws}
The separating hyperplanes constructed in the proof
of Proposition~\ref{th-n2} can be
written as~\eqref{hyperplane},
with $w'\geq w''$ and $d'\geq d''$,
so that,
by the same argument as in Remark~\ref{rk-ineq},
$H$ in~\eqref{hyperplane} may be rewritten as
$H= \bset{v\in \sK^I}\mset{\< w', x\>\oplus d' \leq \< w'', x\>\oplus d''}\eset  \enspace$.
\end{remark}
\begin{proof}
First, we can assume that 
$\support y\subset \support C $, which, by~\eqref{e-sup=},
means that
\begin{align}
\support y= \support C 
\label{e-sup-incl} \enspace .
\end{align}
Otherwise, there is an index
$i\in I$ such that $y_i\neq\zero$ and $v_i=\zero$, for
all $v\in C$, so that the hyperplane of equation
$v_i=\zero$ contains $C$ and not $y$. 

We can also assume that 
\begin{align}
\support C=I
\label{e-a-1} \enspace .
\end{align}
Indeed,  if $\support C \neq I$, we set
$J:=\support C$, and consider the restriction map
$r: \sK^I\to \sK^J$, which sends
a vector $x\in \sK^I$ to $r(x)=(x_j)_{j\in J}$. 
We have $\support r(y)\subset \support r(C)=J$.
Assuming that the theorem is proved when~\eqref{e-a-1}
holds, we get vectors $w',w''\in \sK^J$ 
and scalars $d',d''\in \sK$ such that the affine hyperplane 
$H=\bset x\in \sK^J\mset \<w',x\>\oplus d'= \<w'',x\>\oplus d''\eset$
contains $r(C)$ and not $r(y)$. Let $\hat w'$ and $\hat w''$
denote the vectors obtained by completing $w'$ and $w''$ by zeros.
Then, the hyperplane
$\hat H= \bset x\in \sK^I\mset \<\hat w',x\>\oplus d'= \<\hat w'',x\>\oplus d''\eset$ contains $C$ and not $y$. 

It remains to show Proposition~\ref{th-n2} when the equalities
~\eqref{e-sup-incl}, ~\eqref{e-a-1} hold.
Define the complete convex set $\bar C$ as in~\eqref{e-def-barc}.
It follows from~\eqref{e-fond} that $y\not \in \bar C$.
Therefore, defining $Q_{\bar C}(y)$ and $\nu_{\bar C}(y)$
as in~\eqref{e-formul2}, with $C$ replaced by $\bar C$, we get that the set
$H$ of~\eqref{e-separh-convex}, where $C$ is replaced by $\bar C$,
contains $C$ and not $y$. 
By \eqref{e-sup-incl} and~\eqref{e-a-1},
we have $\support y=I$, so $\ldual{y}\in  \sK^{I}$. 
Also, by~\eqref{e-downbarqc}
and~\eqref{e-fond}, $Q_{\bar C}(y)\in \sK^I$.
Moreover $\nu_{\bar C}(y)\in \sK$ ($\sK$ is conditionally
complete and $\nu_{\bar C}(y)$ is the
sup of a family of elements bounded from above by the unit). 
If 
\begin{align}
\label{e-toc}
(Q_{\bar C}(y))_i\neq\zero,\;\forall i\in I,\;\;
\mrm{and}\;\;
\nu_{\bar C}(y)\neq\zero \enspace ,
\end{align}
we will have $\ldual{Q_{\bar C}(y)}\in \sK^I$,
and $\ldual{\nu_{\bar C}(y)}\in \sK$, 
and the set $H\cap \sK^I$, where $H$ is as in~\eqref{e-separh-convex},
will be an affine hyperplane of $\sK^I$.
In order to show~\eqref{e-toc}, take
any $i\in I$. Since the equalities~\eqref{e-sup-incl} and~\eqref{e-a-1}
hold, we can find $v\in C$ such that $v_i\neq\zero$, and 
thanks to Assumption~(A), $v\lambda \leq y$,
for some $\lambda\in \sKp$. Hence,
\[
\nu_{\bar C}(y)\geq v\lres y\wedge\unit 
\geq v \lres (v\lambda) \wedge \unit \geq \lambda \wedge \unit>\zero
\enspace ,
\]
and 
\[
(Q_{\bar C}(y))_i \geq v_i(v\lres y \wedge \unit)
\geq v_i (\lambda\wedge \unit)
>\zero\enspace,
\]
which shows~\eqref{e-toc}.
\end{proof}
When $C$ is a semimodule,
the condition that $C$ is stable under filtered
infs (which is implied by the condition
that $C$ is order closed) can be dispensed with.
\begin{proposition}\label{th-n1}
Let $C$ be a subsemimodule of $\sK^I$,
and $y\in \sK^I\setminus C$.
Assume that $C$ is stable under directed sups, 
and that Assumption (A) is satisfied.
Then, there is a linear hyperplane of $\sK^I$
which contains $C$ and not $y$. 
\end{proposition}
\begin{proof}
We reproduce the proof of Proposition~\ref{th-n2},
using directly~\eqref{e-linear}, where $V=\bar C$,
and noting that, by Proposition~\ref{prop-cbar2},
\eqref{e-fond} holds as soon
as $C$ is stable under directed sups, when $C$
is a semimodule.
\end{proof}
In Proposition~\ref{th-n2}, we required the convex
set to be order closed, but the separating sets, namely
the affine hyperplanes of $\sK^I$,
where $I$ is infinite, need not be order closed,
as shown by the following counter-example.
\begin{example}
Let $I=\NN$, $\sK=\rmax$, and let us separate $y=(0,1,0,1,0,1\ldots)$
from the convex set $C=\{(0,0,0,\ldots)\}$ using Proposition~\ref{th-n2}.
We obtain the affine
hyperplane $H=\bset x\in \rmax^\NN \mset a(x)\oplus 0
=b(x)\oplus 0\eset$, where $a(x)=x_0\oplus (-1)x_1\oplus x_2\oplus(-1)x_3\oplus\cdots$,
and $b(x)=x_0\oplus x_1\oplus x_2\oplus\cdots$.
Consider the decreasing sequence $y^\ell\in \rmax^\NN$,
such that $y^\ell_{2i+1}=2$, for all $i\in \NN$, 
and $y^\ell_{2i}=1$, for $i\leq \ell$, and $y^\ell_{2i}=2$,
for $i>\ell$, so that $y^\ell\in H$ for
all $\ell$. We have $\inf_\ell y^\ell=y$, where $y_{2i+1}=2$
for all $i\in \NN$, and $y_{2i}=1$, for all $i\in \NN$.
Since $y\not\in H$,
$H$ is not stable under filtered infs. 
\end{example}
Of course, this pathology vanishes
in the finite dimensional case.
\begin{proposition}\label{affineareclosed}
Affine hyperplanes of $\sK^n$ are closed in the order topology.
\end{proposition}
\begin{proof}
This follows readily from Lemma~\ref{lem-birkhoff}.
\end{proof}
The following example shows
that the archimedean assumption is useful in Proposition~\ref{th-n1}.
\begin{example}
Consider the semimodule $C=\{(-\infty,-\infty)\}
\cup \bset (x_1,x_2)\in\RR\times \RR\mset x_1\geq x_2
\eset\subset\rmax^2$, which is stable under directed sups,
and consider the point $y=(0,1)\not\in C$,
with $\support y=\support C=\{1,2\}$, so that Assumption~(A)
is satisfied. The proof of Theorem~\ref{th-n1} allows us to separate $y$
from $C$ by the linear hyperplane:
\begin{align}
H=\bset (x_1,x_2)\in \rmax^2\mset x_1 \oplus (-1)x_2 = x_1 \oplus x_2
\eset \enspace .
\end{align}
However, consider now $y=(0,-\infty)\not\in C$,
which does not satisfy Assumption~(A). We cannot separate
$y$ from $C$ by a linear (or affine) hyperplane,
because such an hyperplane would be closed
in the order (=usual) topology of $\rmax^2$
(by Proposition~\ref{affineareclosed}) whereas
$y$ belongs to the closure of $C$ in this topology.
Thus Assumption~(A) cannot be ommited in Proposition~\ref{th-n2}.
\end{example}
\subsection{Projectors onto closed semimodules of $\protect\mathcal{K}^n$}
In order to show that in the finite dimensional case,
Assumption~(A) is not needed in Proposition~\ref{th-n2},
we establish some continuity property for 
projectors onto closed semimodules of $\sK^n$.

If $V$ is a subsemimodule of $\sK^n$,
we define $\bar V\subset\kbar^n$ as in~\eqref{e-def-barc}
(the condition $\bigoplus_{\ell\in \ELL} \lambda_\ell=\unit$
can be dispensed with, since $V$ is a semimodule),
together with the projector 
$P_\barV: \kbar^n \to \kbar^n,$
\[
P_\barV(x)= \maxx\bset v\in \barV\mset v\leq x\eset\;.
\]
Since $\{v\}_{v\in V}$ is a generating
family of the complete semimodule $\barV$,
it follows from~\eqref{e-proj} that 
\begin{align}
P_\barV(x)= \bigoplus_{v\in V} v(v\lres x) \enspace .
\label{eq-charac-P}
\end{align}
\begin{proposition}
If $V$ is a subsemimodule of $\sK^n$,
that is stable under directed sups,
then the projector $P_\barV$ from $\kbar^n$ onto $\barV$
admits a restriction $P_V$ from $\sK^n$ to $V$.
\end{proposition}
\begin{proof}
If $y\in \sK^n$, $P_\barV(y)\leq y$ 
also belongs to $\sK^n$, so that by~\eqref{e-fond},
$P_\barV(y)\in \barV\cap \sK^n=V$.
\end{proof}
\begin{definition}
We say that a map $f$ from $S$ to an ordered
set $T$ \new{preserves directed sups}
(resp.\ \new{preserves filtered infs})
if $f(\supp D)=\supp f(D)$ (resp.
$f(\inff F)=\inff f(F)$ 
for all directed subsets $D\subset S$
bounded from above (resp.\ for
all filtered subsets $F\subset S$ bounded
from below). 
\end{definition}
\begin{proposition}\label{theo-cont}
If \(V\) is a subsemimodule of \(\sK^n\)
stable under directed sups and filtered infs,
then \(P_V\) preserves directed sups and filtered
infs. 
\end{proposition}
\begin{proof}
Let $F$ denote a filtered subset of $V$,
and $x=\inff F$.
Then, by $\inff F\geq P_V(\inff F)$,
and since $P_V$ is isotone, we have
\begin{align}
\label{e-pv}
P_V(x)&=P_V(\inff F)   \geq  P_V(\inff P_V(F)) \enspace .
\end{align}
Furthermore, since $P_V$ is isotone, $P_V(F)$ is filtered
(indeed, if $F'$ is any finite subset of $P_V(F)$, we 
can write $F'=P_V(F'')$ for some finite subset $F''\subset F$;
since $F$ is filtered, $F''$ has a lower bound $t\in F$,
and since $P_V$ is isotone, $P_V(t)\in P_V(F)$ is a lower bound
of $F'=P_V(F'')$, which shows that $P_V(F)$ is filtered).
Hence, $\inff P_V(F)\in V$ because $V$ is stable under
filtered infs. Since $P_V$ fixes
$V$, $P_V(\inff P_V(F))=\inff P_V(F)$, and
we get from~\eqref{e-pv},
\(
P_V(\inff F)\geq \inff P_V(F)\).
The reverse inequality is an immediate
consequence of the isotony of \(P_V\).

Consider now a directed subset $D\subset V$
bounded from above, and $x=\supp D=\bigoplus_{y\in D}y$. 
We first show that for all $v\in\sK^n$,
\begin{align}
\bigoplus_{y\in D} v(v\lres y) = v(v\lres \bigoplus_{y\in D} y)\enspace .
\label{e-c1}
\end{align}
We shall assume that $v\neq\zero$ (otherwise, the equality
is trivial). Since $D$ is directed,
the net $\{y\}_{y\in D}$ order converges to $\bigoplus_{y\in D} y$,
and by Corollary~\ref{cor-order}, this implies
that $\{v(v\lres y)\}_{y\in D}$ order converges
to $v(v\lres \bigoplus_{y\in  D} y)$. Since $y\mapsto
v(v\lres y)$ is isotone, $\{v(v\lres y)\}_{y\in D}$ order converges
to its sup. (Indeed, let $\varphi$
denote an isotone map from $D$ to a conditionally complete ordered
set, such that $\{\varphi(y)\}_{y\in D}$
is bounded from above, and let us show more generally that
$\{\varphi(y)\}_{y\in D}$ order converges to its sup. Observe that 
$\sup_{z\geq y}\varphi(z)=\supp \varphi(D)$ is independent of $y\in D$
because $D$ is directed and $\varphi$ is isotone.
Then, $\limsup_{y\in D}\varphi(y)
=\inf_{y\in D} \sup_{y'\geq y}\varphi(y')=\supp \varphi(D)$.
Also, since $\varphi$ is isotone,
$\liminf_{y\in D}\varphi(y)=\sup_{y\in D}\inf_{y'\geq y}\varphi(y')
=\sup_{y\in D}\varphi(y)=\supp\varphi(D)$,
which shows that $\{\varphi(y)\}_{y\in D}$ order converges to its sup.)
So, \eqref{e-c1} is proved.

Using~\eqref{eq-charac-P}, we get
\begin{align*}
\bigoplus_{y\in D} P_V(y)=\bigoplus_{y\in D} \bigoplus_{v\in V}
v(v\lres y)  = \bigoplus_{v\in V}\bigoplus_{y\in D}v(v\lres y)  
&=\bigoplus_{v\in V}v\big(v\lres (\bigoplus_{y\in D} y))\big) \quad\mrm{(by~\eqref{e-c1})}\\
&=P_V(x) \enspace .
\end{align*}
\end{proof}
\begin{remark}
Proposition~\ref{theo-cont} does not extend
to semimodules of the form $\sK^I$, where
$I$ is an infinite set. Indeed, take $I=\NN$,
$\sK=\rmax$, and let $V$ denote the semimodule spanned
by the vector $v=(0,0,0,\ldots)$.
For all $x=(x_0,x_1,x_2,\ldots)\in\rmax^{\NN}$,
we have $P_V(x)=(\lambda(x),\lambda(x),\ldots)$,
where $\lambda(x):=\inff_{i\in\NN}x_i$. Consider 
now the sequence $y^k\in \rmax^{\NN}$, such that
$y^k_i=0$ if $k\leq i$, and $y^k_i=-1$, otherwise.
Then, $\{y^k\}_{k\in \NN}$ is a non-decreasing sequence
with supremum $v$. We have $\lambda(y^k)=-1$, but $\lambda(v)=0$,
which shows that $P_V$ does not preserve directed sups.
\end{remark}
The proof of Theorem~\ref{th-1} will rely on the
following corollary of Proposition~\ref{theo-cont}.
\begin{corollary}\label{cor-key}
If $V$ is a subsemimodule of $\sK^n$ stable under directed
sups and filtered infs, 
and if $y\in\sK^n\setminus V$, then, there is a vector $z\geq y$
with coordinates in $\sK\setminus\{\zero\}$,
such that:
\begin{align*}
y\not\leq P_V(z) \enspace .
\end{align*}
\end{corollary}
\begin{proof}
Let $Z$ denote the set of vectors
$z\geq y$ with coordinates in $\sK\setminus\{\zero\}$.
Let us assume by contradiction that
\begin{align}
y\leq P_V(z) ,\quad \forall z\in Z \enspace .\label{e-leqf}
\end{align}
Since by Proposition~\ref{theo-cont}, $P_V$ preserves
filtered infs, we get from~\eqref{e-leqf}:
\[
y\leq \inff_{z\in Z} P_V(z)= P_V(\inff_{z\in Z} z)
=P_V(y)
\enspace. 
\]
Since $y\geq P_V(y)$ holds trivially, 
$y=P_V(y)$, hence, $y\in V$, a contradiction.
\end{proof}
\subsection{Separation theorem for closed convex subsets of $\protect\mathcal{K}^n$}
The following finite dimensional separation theorem
extends an earlier result of Zimmermann~\cite{zimmerman77}.
Recall that when $\sK=\rmax$, the order topology on $\sK^n$
is the usual topology on $\rmax^n$.
\begin{theorem}\label{th-1}
Let $C$ denote a convex subset of $\sK^n$
that is closed for the order topology of $\sK^n$,
and let $y\not\in C$. Then, there exists
an affine hyperplane containing $C$ and not $x$.
\end{theorem}
We shall need the following lemma:
\begin{lemma}\label{support=}
If $C$ is a semimodule, and
if $\support y=\support C$, then $\support P_C(y)=\support y$.
\end{lemma}
\begin{proof}
Since $P_C(y)\leq y$, $\support P_C(y)\subset \support y$.
Conversely, pick any $i\in \support y$. 
Since $\support y=\support C$, we can find
$v\in C$ such that $\support v\subset \support y$ and $i\in\support v$.
Then, by~\eqref{e-def-lres}, 
$v\lres y=\inff_{j\in \support v}v_j^{-1}y_j \neq\zero$,
and since $P_C(y)\geq v(v\lres y)$, $(P_C(y))_i\neq\zero$,
which shows that $\support y\subset \support P_C(y)$.
\end{proof}
We showed in the first part of the proof of Proposition~\ref{th-n2}
that we can always assume that
\begin{align}
\support y \subset \support C =\{1,\ldots,n\} \enspace .\label{inclusion}
\end{align}
We next prove Theorem~\ref{th-1} in the special
case where $C$ is a semimodule, and 
then, we shall derive Theorem~\ref{th-1}, in general.
\begin{proof}[Proof of Theorem~\ref{th-1} when $C$ is a semimodule]
The proof relies on a perturbation argument.
Pick a vector $z\geq y$ with coordinates in $\sK\setminus\{\zero\}$,
(hence, by~\eqref{inclusion}, $\support z=\{1,\ldots,n\}=\support C$),
and define, as in~\eqref{e-def-seph},\eqref{e-linear}:
\begin{align*}
H(z) = 
\bset x\in \sK^n \mset x\lres z  = x\lres P_C(z) \eset 
=
\bset x\in \sK^n \mset \< \ldual{z}, x\>  =
\<\ldual{P_C(z)}, x\>  \eset \enspace .
\end{align*}
We will show that $H(z)$ is a (linear)
hyperplane, and that one can choose the above $z$ so that
$H(z)$ contains $C$ and not $y$.

It follows from 
$\support z=\support C$ and Lemma~\ref{support=}
that $\support P_C(z)=\support z=\{1,\ldots,n\}$.
Since $\ldual{u}\in \sK^n$ for all vectors
$u$ of $\sK^n$ with coordinates different from $\zero$,
$\ldual{z}$ and $\ldual{(P_C(z))}$
belong to $\sK^n$, which shows that $H(z)$
is an hyperplane.

By Theorem~\ref{th-separ}, $H(z)$ contains $C$.
Let us check that:
\begin{align}
\label{e-imp}
x\in H(z) \implies P_C(z) \geq x \enspace .
\end{align}
Recall the classical residuation
identity
\begin{align}
x(x\lres x)=x
\label{e-res}
\end{align}
(this can be shown by applying the first identity in~\eqref{ffshf=f}
to the map $f:\sK\to \sK^n, f(\lambda)=x\lambda$).
If $x\in H(z)$, we have $x \lres  z = x\lres P_C(z)$,
and by~\eqref{e-1},
$P_C(z) \geq x(x\lres z)$.
Using $z\geq x$ and~\eqref{e-res},
we get 
$P_C(z)\geq x(x\lres z)\geq x(x\lres x)=x$,
which shows~\eqref{e-imp}.

By Corollary~\ref{cor-key}, 
there is a vector $z\geq y$ with
entries in $\sKp$ such that $y\not\leq P_C(z)$,
and by~\eqref{e-imp},
the hyperplane $H(z)$ associated to this $z$ contains
$C$ and not $y$.
\end{proof}
Associate to a convex set $C\subset \sK^n$ the semimodule:
\[
V_C:= \bset (x\lambda,\lambda) \mset x\in C,\;\lambda\in \sK \eset
\subset \sK^{n+1}
\enspace .
\]
We denote by $\clo(V_C)$ the closure of $V_C$
for the order topology of $\sK^{n+1}$.
We shall need the following:
\begin{lemma} If $C$ is a convex subset of $\sK^n$
closed for the order topology, 
\begin{align}
\clo(V_C) \subset V_C \cup (\sK^n\times \{\zero\})  \enspace .
\label{e-closure}
\end{align}
\end{lemma}
\begin{proof}
It suffices to show that $V_C\cup (\sK^n\times \{\zero\})$ is 
closed in $\sK^{n+1}$ for the order topology. Take a net 
$\{(z_\ell,\lambda_\ell)\}_{\ell\in \ELL}\subset V_C\cup (\sK^n\times \{\zero\})$,
with $z_\ell\in \sK^n,\lambda_\ell\in\sK$, order converging to 
some $(z,\lambda)$, with $z\in \sK^n,\lambda\in\sK$.
We only need to show that if $\lambda\neq\zero$,
$(z,\lambda)\in V_C$. Since $\lambda\neq\zero$, 
replacing $\ELL$ by a set of the form $\bset \ell\in \ELL\mset \ell \geq \ell_0\eset$,
we may assume that $\lambda_\ell\neq \zero$, 
for all $\ell\in \ELL$. Then, $(z_\ell,\lambda_\ell)\in V_C$,
which implies that $z_\ell\lambda_\ell^{-1}\in C$. Since
$z_\ell$ order converges to $z$, and $\lambda_\ell^{-1}$
order converges to $\lambda^{-1}$, by Lemma~\ref{lem-birkhoff},
$z_\ell\lambda_\ell^{-1}$ order-converges to $z\lambda^{-1}$.
Since, by our assumption, $C$ is closed for
the order topology, $z\lambda^{-1}\in C$, which shows
that $(z,\lambda)\in V_C$. So, $V_C\cup  (\sK^n\times \{\zero\})$ is 
closed for the order topology.
\end{proof}
\begin{proof}[Derivation of the general case of Theorem~\ref{th-1}]
Let us take $y\in\sK^n\setminus C$.
We note that 
by~\eqref{e-closure},
$(y,\unit)\not\in \clo(V_C)$.
Applying Theorem~\ref{th-1}, which is already proved
in the case of closed semimodules, to $\clo(V_C)$,
which is a semimodule thanks to Lemma~\ref{l-1},
we get a linear hyperplane
$H=\bset \bar x\in \sK^{n+1}\mset \<w',\bar x\> =\<w'',\bar x\>\eset$,
where $w',w''\in \sK^{n+1}$, such that 
\begin{align}\label{separc}
x\in C\implies (x,\unit)\in H,\qquad \text{ and } (y,\unit)\not\in H
\enspace .
\end{align}
Introducing $z'=(w'_i)_{1\leq i\leq n}$,
and $z''=(w''_i)_{1\leq i\leq n}$, we see
from~\eqref{separc} that the affine
hyperplane:
\[
\bset{x\in \sK^n}\mset{\<z',x\>\oplus w'_{n+1}= \<z'',x\>\oplus w''_{n+1}}\eset
\]
contains $C$ and not $y$.
\end{proof}
\begin{remark}
We needed to introduce the closure $\clo(V_C)$
in the proof of Theorem~\ref{th-1}
because $V_C$ need not be closed
when $C$ is closed and convex.
Indeed, consider $C=\bset{x\in \rmax}\mset{x\geq 0}\eset$.
We have
$V_C=\{(\zero,\zero)\}\cup\bset (x\lambda,\lambda)\mset x\geq 0,\lambda\in \RR \eset$ and $\clo(V_C)=(\rmax \times \{\zero\})
\cup \bset (x\lambda,\lambda)\mset x\geq 0,\lambda\in \RR \eset$.
\end{remark}
\begin{example}\label{ex-ncex}
When applied to Example~\ref{cex},  the proof of Theorem~\ref{th-1}
shows that for $k$ large enough, the hyperplane $H$
in~\eqref{e-sep-h} separates the point $(0,-\infty)$ from 
the convex set of Figure~\ref{fig-todo}. 
The method of~\cite{shpiz}, which requires
that the vector to separate from a convex set should have
invertible entries in order to apply a normalization argument,
does not apply to this case.
\end{example}
\section{Convex functions over idempotent semifields}\label{sec-fconvex}
We say that a map $f:\sK^n\to\kbar$ is \new{convex}
if its epigraph is convex. 
By \cite[Theorem 1]{zimmer}, 
$f$ is convex if, and only if,
\[
(x,y\in \sK^n,\alpha ,\beta \in \sK,\alpha \oplus \beta
=\unit)\Rightarrow f(x\alpha\oplus y\beta)\leq f(x)\alpha \oplus f(y)\beta
\enspace .
\]
Additionally, by \cite[Theorem 2]{zimmer}, 
the (lower) level sets 
\[
S_t(f)=\bset x\in X\mset f(x)\leq t\eset\qquad \quad (t\in \kbar)
\]
of $f$ are convex subsets of $\sK^n$. 
When $\sK=\rmax$, we say that $f$ is \new{max-plus convex}.
Convex functions may of course be defined
from an arbitrary $\sK$-semimodule $X$
to $\kbar$: we limit our attention to $X=\sK^n$
since the proof of the main result below relies on 
Theorem~\ref{th-1} which is stated for $\sK^n$.

The following immediate proposition shows that
the set of convex functions is a 
complete subsemimodule of the complete semimodule
of functions $\sK^n\to \kbar$:
\begin{proposition}
\label{lsupmax}The set of all convex functions is stable 
under (arbitrary) pointwise sup, and under multiplication
by a scalar (in $\kbar$).\qed 
\end{proposition}
We defined in the introduction
\new{$U$-convex} functions
and sets,
when $U\subset \rbar^X$,
see Equations~\eqref{uconv} and~\eqref{uconv2}.
When more generally $U\subset \kbar^X$,
we still define \new{$U$-convex functions} by~\eqref{uconv},
and extend~\eqref{uconv2} by saying
that a subset $C\subset \sK^n$ is \new{$U$-convex}
if for all $y\in X\setminus C$, we can find a map $u\in U$
such that 
\begin{align}
u(y) \not\leq \sup_{x\in C} u(x) \enspace .\label{uconv3}
\end{align}
\begin{proposition}\label{prop-inter}
A subset $C\subset X$ is $U$-convex if, and only if,
it is an intersection  of (lower) level sets of maps
in $U$.
\end{proposition}
\begin{proof}
Assume that $C$ is an intersection of
(lower) level sets of maps in $U$, that is,
$C=\bigcap_{\ell\in \ELL} S_{t_\ell}(u_\ell)$,
where $\{u_\ell\}_{\ell\in \ELL}\subset U$, $\{t_\ell\}_{\ell\in \ELL}\subset \kbar$,
and $\ELL$ is a possibly infinite set. 
If $y\in X\setminus C$, $y\not\in S_{t_\ell}(u_\ell)$,
for some $\ell\in \ELL$,
so that $u_\ell(y)\not\leq t_\ell$.
Since $C\subset S_{t_\ell}(u_\ell)$,
we have
$\sup_{x\in C} u_\ell(x)\leq t_\ell$.
We deduce that $u_\ell(y)\not\leq \sup_{x\in C} u_\ell(x)$. 
Hence, $C$ is $U$-convex. 

Conversely, assume that $C$ is $U$-convex,
and let $C'$ denote the intersection
of the sets $S_t(u)$, with $u\in U$ and $t\in\kbar$,
in which $C$ is contained. Trivially, $C\subset C'$. 
If $y\in X\setminus C$, we can find $u\in U$
satisfying~\eqref{uconv3}. Let $t:=\sup_{x\in C}u(x)$.
Then, $C\subset S_t(u)$, and $y\not\in S_t(u)$,
so that $y\not\in C'$.
This shows that $X\setminus C\subset X\setminus C'$.
Thus, $C=C'$.
\end{proof}
The set \(U\) of elementary functions which will prove relevant for our convex functions is the following. 
\begin{definition}
We say that $u:\sK^n\to \sK$
is \new{affine} if $u(x)=\proaff{w}{d}{x}$, for
some $w\in\sK^n$ and $d\in\sK$. We say that $u$
is a \new{difference of affine functions}
if 
\begin{equation}
u(x)= \afff{w}{d}{x}\;,
\label{wefing3}
\end{equation} 
where \(w', w''\) belong to \(\sK^{n}\) and \(d',d''\) to \(\sK\).
\end{definition}

We illustrate in Table~\ref{tab:1} below, and in Figure~\ref{fig:affine}
of \S\ref{sec-intro}, the various shapes
taken by differences of affine functions, when $n=1$, and $\sK=\rmax$.
For simplicity, a generic function in this class is denoted
\begin{displaymath}
y=(a x \oplus b)\ominus (c x \oplus d)
\end{displaymath}
with \(a,b,c,d\in\sK\).
Table~\ref{tab:1} enumerates the four possible situations according to the comparisons of \(a\) with \(c\) and \(b\) with \(d\).
\begin{table}[hbtp]
\begin{displaymath}\begin{array}{c|c|c}
&a\leq c&a > c\\
\hline b\leq d&y=\zero, \quad\forall x& y=\begin{cases}
\zero& \text{if\ }x\leq a\lres d, \\
a x & \text{otherwise}.
\end{cases}\\
\hline b> d&y=\begin{cases}
b& \text{if\ } x<b\lres c, \\
\zero & \text{otherwise}.
\end{cases}&y = a x \oplus b
\end{array}
\end{displaymath}
\caption{The four generic differences of affine functions over \(\rmax\)}\label{tab:1}
\end{table}
Figure~\ref{fig:affine} shows the corresponding plots. 
\begin{proposition}\label{prop-ls}
The (lower) level sets of differences of affine functions are precisely the 
affine hyperplanes of the form:
\begin{align}
\bset x\in \sK^n \mset
\proaff{w''}{d''}{x}  \geq \proaff{w'}{d'}{x}\eset
\enspace,
\label{e-ineq}
\end{align}
where \(w', w''\) belong to \(\sK^{n}\) and \(d',d''\) to \(\sK\).
\end{proposition}
\begin{proof}
This is an immediate consequence of~\eqref{e-1dual}.
\end{proof}
\begin{remark}
The inequality in~\eqref{e-ineq} is equivalent 
to the equality 
$\proaff{w''}{d''}{x}  = \proaff{w'\oplus w''}{d'\oplus d''}{x}$,
which justifies the term ``affine hyperplane''.
\end{remark}
We shall say that $f:\sK^n\to\kbar$ is \new{lower semi-continuous}
if all (lower) level sets of $f$ are closed
in the order topology of $\sK^n$.
\begin{proposition}
\label{theor}
Every difference of affine functions is convex and lower semi-continuous.
\end{proposition}
\begin{proof}
If $u$ is a difference of affine functions,
by Proposition~\ref{prop-ls},
the (lower) level sets of $u$, are affine hyperplanes,
which are closed by Proposition~\ref{affineareclosed}, 
so $u$ is lower semi-continuous.

As mentioned earlier, the function $u(\cdot)=\afff{w}{d}{\cdot}$ is  convex if and only if its epigraph is  convex.
So, we consider two points \((x_1,\lambda_{1})\) and \((x_2,\lambda_{2})\) in the epigraph of \(u\), namely,
\begin{align*}
\lambda_{1}&\geq \afff{w}{d}{x_1}\;,\\
\lambda_{2}&\geq \afff{w}{d}{x_2}\;,
\end{align*}
which, by~\eqref{e-1dual}, is equivalent to
\begin{align*}
\lambda_{1}\oplus\proaff{w''}{d''}{x_1}  &\geq \proaff{w'}{d'}{x_1} \;,\\
\lambda_{2}\oplus\proaff{w''}{d''}{x_2} &\geq \proaff{w'}{d'}{x_2} \;.
\end{align*}
Let \(\alpha\) and \(\beta\) in \(\sK\) be such that \(\alpha\oplus\beta=\unit\). From the previous inequalities, we derive
\begin{displaymath}
\lambda_{1}\alpha\oplus \lambda_{2}\beta\oplus\proaff{w''}{d''}{x_1\alpha\oplus x_2\beta}\geq \proaff{w'}{d'}{x_1\alpha \oplus x_2\beta}\;,
\end{displaymath}
which, by~\eqref{e-1dual}, is equivalent to 
\begin{displaymath}
\lambda_{1}\alpha \oplus  \lambda_{2}\beta\geq u( x_1\alpha \oplus x_2 \beta)\;.
\end{displaymath}
We have proved that \(( x_1\alpha\oplus x_2\beta, \lambda_{1}\alpha\oplus \lambda_{2}\beta)\) belongs to the epigraph of \(u\).
Thus, \(u\) is  convex.
\end{proof}
\begin{corollary}\label{cor-Uconvex}
Let $C\subset \sK^n$. 
The following assertions are equivalent:
\begin{enumerate}
\item
$C$ is a convex subset of $\sK^n$, and it is closed in the order topology;
\item
$C$ is $U$-convex,
where $U$ denotes the set of differences of affine functions $\sK^n\to \sK$,
defined by~\eqref{wefing3}.
\end{enumerate}
\end{corollary}
\begin{proof}
If $C$ is convex and closed, by Theorem~\ref{th-1},
for all $y\in \sK^n\setminus C$,
we can find an hyperplane~\eqref{hyperplane}
containing $C$ and not $y$.
By Remark~\ref{rk-wpgeqws},
we can choose $w',w'',d',d''$ in~\eqref{hyperplane}
so that $w'\geq w''$ and $d'\geq d''$.
Since $\proaff{w'}{d'}{x}=\proaff{w''}{d''}{x}$,
for all $x\in C$, $u(x)=\afff{w}{d}{x}=\zero$
for all $x\in C$, so that $\sup_{x\in C}u(x)=\zero$. 
Since $\proaff{w'}{d'}{y}\neq \proaff{w''}{d''}{y}$,
and $\proaff{w'}{d'}{y}\geq \proaff{w''}{d''}{y}$
because $w'\geq w''$ and $d'\geq d''$, we must
have $\proaff{w'}{d'}{y}\not\leq \proaff{w''}{d''}{y}$.
Then, $u(y)=\afff{w}{d}{y}\not\leq \zero$,
which shows that $C$ is $U$-convex.

Conversely, if $C$ is $U$-convex,
Proposition~\ref{prop-inter}
shows that $C$ is an intersection
of (lower) level sets of differences of affine functions.
By Proposition~\ref{prop-ls}, these (lower) level 
sets all are affine hyperplanes, and a fortiori,
are convex sets.  Moreover, 
by Proposition~\ref{affineareclosed},
affine hyperplanes are closed,
so $C$ is closed and convex.
\end{proof}
\begin{theorem}
\label{tconverse}
A function $f:\sK^n\rightarrow \kbar$
is convex and lower semi-continuous if, and only if,
it is a sup of differences of affine functions,
i.e., a $U$-convex function, where $U$ is the 
set of functions of the form~\eqref{wefing3}.
\end{theorem}
The proof relies on the following extension to the case
of functions with values in a partially ordered set, of
a well known characterization of abstract convexity of functions, in terms of
``separation'' (see \cite[Prop.~1.6i]{dk},  or
\cite[Th.~3.1, Eqn~(3.31)]{ACA}).
\begin{lemma}
For any set $X$, 
and $U\subset \sK^X$, 
a map $f:X\to \kbar$ is $U$-convex if, and only if, for each
$(x,\nu )\in X\times \sK$ such that \(f(x)\not\leq \nu\), there exists $u\in U$ such that 
\begin{equation}
u\leq f,\qquad u(x)\not\leq \nu \enspace.  \label{dol-kur} 
\end{equation}
\end{lemma}
\begin{proof}
Let $g=\sup_{u\in U,\;u\leq f} u$.
Note first that $f$ is $U$-convex if, and only if, 
for all $x \in X$,
$g(x)= f(x)$,
or equivalently,
$g(x)\geq  f(x)$
(the other inequality always holds).
Recall that for all $t\in \kbar$, $\Up(t)=\bset s\in \kbar\mset s\geq t\eset$
denotes the upper set generated by $t$. Trivially: $s\geq t\Leftrightarrow 
\Up(s)\subset \Up(t)\Leftrightarrow \kbar\setminus \Up(s)\supset\kbar
\setminus \Up(t)$. Applying this to $s=g(x)$
and $t=f(x)$, we rewrite $g(x)\geq f(x)$
as 
\begin{align}
f(x)\not\leq \nu \implies g(x)=\sup_{u\in U,\; u\leq f} u(x) \not\leq \nu
\enspace .\label{e-impl}
\end{align}
Since $ \sup_{u\in U,\; u\leq f} u(x) \not\leq \nu$ if, and only
if, $u(x)\not\leq \nu$ for some $u\in U$ such that $u\leq f$,
and since it is enough to check the implication~\eqref{e-impl}
when $\nu\in \sK$ (if $\nu$ is the top element of $\kbar$, 
the implication~\eqref{e-impl} trivially holds), 
the lemma is proved.
\end{proof}
\begin{proof}[Proof of Theorem~\ref{tconverse}]
$2 \Rightarrow 1$. By Proposition~\ref{theor} and Proposition~\ref{lsupmax},
every sup of functions belonging to $U$ is  convex and lower
semi-continuous.

$1\Rightarrow 2$. Assume that $f:\sK^n\rightarrow
\kbar$ is  convex and lower semi-continuous, and let us prove
that $f$ is $U$-convex. 

As mentioned above,  the epigraph of $f$, $\epi f$,
is a  convex closed subset of $\sK^n\times \sK$.
Consider  $(y,\nu )\in \sK^n\times \sK\setminus \epi f$,
so that \(f(y)\not\leq\nu\).
By Theorem~\ref{th-1}, there exist
\((w',\mu',d')\) and \((w'',\mu'',d'')\) in \(\sK^{n}\times \sK\times \sK\)
with $(w',\mu',d')\geq (w'',\mu'', d'')$
such that 
\[
\proaff{(w',\mu')}{d'}{(z,\lambda )}\leq \proaff{(w'',\mu'')}{d''}{(z,\lambda )}\;,\quad\forall(z,\lambda )\in \epi f\;,
\]
\[
\proaff{(w',\mu')}{d'}{(y,\nu)} \not\leq\proaff{(w'',\mu')}{d''}{(y,\nu)}\;,  
\]
that is, 
\begin{equation}
\proaff{w'}{\mu'\lambda\oplus d'}{z}\leq \proaff{w''}{\mu''\lambda\oplus d''}{z}\;,\quad\forall(z,\lambda )\in \epi f\;,
\label{prima1}
\end{equation}
\[
\proaff{w'}{\mu'\nu\oplus d'}{y} 
\not\leq \proaff{w''}{\mu''\nu\oplus d''}{y}\;.
\]
Since the function identically equal to the top element of $\kbar$
is trivially $U$-convex, we shall assume that 
$f\not\equiv \maxx\kbar$, i.e., $\epi f\neq\emptyset$. Then, we claim that 
\begin{equation}
\mu'=\mu''\;.  \label{claimut}
\end{equation}
Indeed, since $\mu '\geq \mu''$, we may assume that $\mu'\neq\zero$.
Then, taking $(z,\lambda )\in \epi f$, 
with $\lambda $ so large that $\langle z,w'\rangle \oplus \lambda
\mu'\oplus d'=\lambda \mu'$, 
from \eqref{prima1} we
obtain $\lambda \mu '\leq 
\langle z,w''\rangle \oplus \lambda \mu''\oplus d''$.
Then, we cannot have $\mu''=\zero$
(otherwise, $\lambda \mu'$ would be bounded
above independently of $\lambda$). 
Therefore, for $\lambda$ large enough
$\langle z,w''\rangle \oplus \lambda \mu''\oplus d''=\lambda\mu''$,
hence, $\lambda\mu'\leq \lambda\mu''$, which,
by $\mu'\geq\mu''$, implies 
$\lambda\mu'=\lambda \mu''$, and multiplying by $\lambda^{-1}$,
we get~\eqref{claimut}.

Hence, by \eqref{prima1}--\eqref{claimut}, we have
\begin{equation}
\proaff{w'}{\mu'\lambda\oplus d'}{z}\leq \proaff{w''}{\mu'\lambda\oplus d''}{z}\;,\quad\forall(z,\lambda )\in \epi f\;,   \label{prima2}
\end{equation}
\begin{equation}
\proaff{w'}{\mu'\nu\oplus d'}{y} \not\leq
\proaff{w''}{\mu'\nu\oplus d''}{y}\;. \label{doua2}
\end{equation}

Let 
\[\dom f=\bset x\in \sK^n\mset f(x) \neq\top\kbar\eset
=\bset x\in \sK^n \mset (x,\lambda)\in \epi f \;\mrm{for}\;\mrm{some}\;\lambda\in\sK\eset
\enspace .
\]
We claim that if $y\in \dom f$, then $\mu'\neq\zero$.
Indeed, if $\mu'=\zero$, then \eqref{prima2} and \eqref{doua2} become 
\begin{equation}
\proaff{w'}{d'}{z}\leq \proaff{w''}{d''}{z}\;,\quad\forall(z,\lambda )\in \epi f\;,  \label{cuc1}
\end{equation}
\begin{equation}
\proaff{w'}{d'}{y} \not\leq \proaff{w''}{d''}{y}\;;  \label{cuc2}
\end{equation}
but, \eqref{cuc1}, together with $(w',d')\geq (w'',d''),$ yields $\proaff{w'}{d'}{z}=\proaff{w''}{d''}{z}$, for all $(z,\lambda )\in \epi f$, that is, $\proaff{w'}{d'}{\cdot}=\proaff{w''}{d''}{\cdot}$
when restricted to arguments lying in \(\dom f\), in contradiction with \eqref{cuc2}, provided \(y\in \dom f\).
This proves the claim $\mu' \neq \zero $ in this case.

In \eqref{prima2}, \eqref{doua2}, we
may now assume that $\mu'=\unit$. Indeed, multiply by $(\mu ')^{-1}$,
and rename $(\mu')^{-1} w',(\mu ')^{-1} w^{'' }, (\mu ')^{-1} d',(\mu')^{-1} d''$ as $w',w'',d',d''$
respectively.  Now \eqref{prima2} and \eqref{doua2} read
\begin{equation}
\proaff{w'}{\lambda\oplus d'}{z}\leq \proaff{w''}{\lambda\oplus d''}{z}\;,\quad \forall (z,\lambda )\in \epi f\;,  \label{prima3}
\end{equation}
\begin{equation}
\proaff{w'}{\nu\oplus d'}{y} \not\leq\proaff{w''}{\nu\oplus d'}{y}\;.
\label{doua3}
\end{equation}

Equation \eqref{prima3} implies that
\[
\proaff{w'}{d'}{z}\leq \proaff{w''}{\lambda\oplus d''}{z}\;,\quad \forall (z,\lambda )\in \epi f\;,
\]
whence, by \eqref{e-1dual}, 
\begin{equation}
\afff{w}{d}{z}\leq \lambda\;,\quad
\forall (z,\lambda )\in \epi f\;,  \label{bunbun}
\end{equation}
and hence, defining 
\begin{equation}
u:=\afff{w}{d}{\cdot}\in U\;,  \label{defu}
\end{equation}
and using that $f(z)=\minn\bset\lambda\mset (z,\lambda )\in \epi f\eset$,
from \eqref{bunbun} we see that \(f(z)\geq u(z)\) for \(z\in\dom f\); for \(z\notin \dom f\), \(f(z)=\top\kbar\) and this inequality is trivial,
thus we have obtained the first half of \eqref{dol-kur}.

From~\eqref{doua3}, we deduce that
\[
\proaff{w'}{d'}{y} \not\leq \proaff{w''}{d''\oplus \nu}{y}\
\]
(because $a\leq b\oplus \nu \implies a\oplus \nu \leq b\oplus \nu$),
whence, by \eqref{e-1dual}
(in fact, its equivalent negative form) and \eqref{defu}, we obtain 
\[
u(y)=\afff{w}{d}{y} \not\leq \nu , 
\]
that is, the second part of \eqref{dol-kur}.

For the proof to be complete, we have to handle the case when \(y\notin \dom f\) which implies that \((y,\nu)\notin \epi f\).
The previous arguments hold true up to a certain point
when we cannot claim that \(\mu'\neq\zero\).
Either $\mu'\neq\zero$ indeed, and the proof is completed as previously, 
or $\mu'=\zero$, and then \eqref{prima3}--\eqref{doua3} boil down to 
\begin{equation}
\proaff{w'}{d'}{z}\leq \proaff{w''}{d''}{z}\;, \quad\forall z\in \dom f\;,  \label{prima4}
\end{equation}
\begin{equation}
\proaff{w'}{d'}{y} \not\leq \proaff{w''}{d''}{z},
\label{doua4}
\end{equation}
without having to redefine the original \((w',d',w'',d'')\). 
For any \(\alpha\in\sK\setminus\{\zero\}\), define the functions 
\begin{align}
u_{\alpha}(\cdot)&=\alpha(\proaff{w'}{d'}{\cdot})\ominus \alpha(\proaff{w''}{d''}{\cdot})=\alpha\big(\afff{w}{d}{\cdot}\big)\; ,  
\label{e-mult}
\end{align}
which all belong to \(U\). Because of \eqref{prima4}, all those functions are identically equal to \(\zero\) over \(\dom f\), hence they are trivially less than or equal to \(f\) over \(\dom f\) but also over the whole \(\sK^{n}\).
On the other hand, because of \eqref{doua4}, 
\begin{displaymath}
\afff{w}{d}{y}>\zero\;.
\end{displaymath}
Multiplying this strict inequality by a large enough
$\alpha$, and using~\eqref{e-mult},
we see that given \(\nu\), there exists \(\alpha\) for which
\(u_{\alpha}(y)>\nu\), and a fortiori, $u_{\alpha}(y)\not\leq\nu$.
The proof is now complete.
\end{proof}
\begin{remark}
By Remark~\ref{rk-wpgeqws}, 
Theorem~\ref{tconverse} remains
valid if the set $U$ is replaced by the subset
of the functions in~\eqref{wefing3} such that 
$w'\geq w''$ and $d'\geq d''$.
\end{remark}
An illustration of Theorem~\ref{tconverse}
has been given in Figure~\ref{fig:fconvex}
in the introduction: the figure shows
a convex function over \(\rmax\), 
together with its supporting hyperplanes (which are epigraphs
of differences of affine functions, whose shapes
were already shown in Figure~\ref{fig:affine}).

\newcommand{\etalchar}[1]{$^{#1}$}


\begin{thebibliography}{GHK{\etalchar{+}}80}
\providecommand{\url}[1]{\texttt{#1}}
\providecommand{\urlprefix}{URL }
  \providecommand{\doi}[1]{Eprint \href{http://dx.doi.org/#1}{doi:#1}}
\providecommand{\arxiv}[2][]{Also \href{http://www.arXiv.org/abs/#2}{arXiv:#2}}

\bibitem[AS03]{akiansinger}
M.~Akian and I.~Singer.
\newblock Topologies on lattice ordered groups, separation from closed downward
  sets, and conjugations of type {Lau}.
\newblock \emph{Optimization}, 52(6):629--672, 2003.

\bibitem[BCOQ92]{bcoq}
F.~Baccelli, G.~Cohen, G.~Olsder, and J.~Quadrat.
\newblock \emph{Synchronization and Linearity --- an Algebra for Discrete Event
  Systems}.
\newblock Wiley, 1992.

\bibitem[Bir67]{birkhoff40}
G.~Birkhoff.
\newblock \emph{Lattice Theory}, volume XXV of \emph{American Mathematical
  Society Colloquium Publications}.
\newblock A.M.S, Providence, Rhode Island, 1967.
\newblock (third edition).

\bibitem[BJ72]{blyth72}
T.~Blyth and M.~Janowitz.
\newblock \emph{Residuation Theory}.
\newblock Pergamon press, 1972.

\bibitem[CG79]{cuning}
R.~Cuninghame-Green.
\newblock \emph{Minimax Algebra}.
\newblock Number 166 in Lecture notes in Economics and Mathematical Systems.
  Springer, 1979.

\bibitem[CGQ96]{CGQ96a}
G.~Cohen, S.~Gaubert, and J.~Quadrat.
\newblock Kernels, images and projections in dioids.
\newblock In \emph{Proceedings of WODES'96}. IEE, Edinburgh, August 1996.

\bibitem[CGQ97]{CGQ97}
G.~Cohen, S.~Gaubert, and J.~Quadrat.
\newblock Linear projectors in the max-plus algebra.
\newblock In \emph{5th IEEE Mediterranean Conference on Control and Systems}.
  Paphos, Cyprus, 1997.

\bibitem[CGQ99]{ccggq99}
G.~Cohen, S.~Gaubert, and J.~Quadrat.
\newblock Max-plus algebra and system theory: where we are and where to go now.
\newblock \emph{Annual Reviews in Control}, 23:207--219, 1999.
\newblock \doi{10.1016/S1367-5788(99)90091-3}.

\bibitem[CGQ01]{praha}
G.~Cohen, S.~Gaubert, and J.~Quadrat.
\newblock Duality of idempotent semimodules.
\newblock In \emph{Proceedings of the Workshop on Max-Plus Algebras, IFAC
  SSSC'01}. Elsevier, Praha, 2001.

\bibitem[CGQ04]{cgq02}
G.~Cohen, S.~Gaubert, and J.~Quadrat.
\newblock Duality and separation theorem in idempotent semimodules.
\newblock \emph{Linear Algebra and Appl.}, 379:395--422, 2004.
\newblock \doi{10.1016/j.laa.2003.08.010}.
\newblock \arxiv{math.FA/0212294}.

\bibitem[CKR84]{cao84}
Z.~Cao, K.~Kim, and F.~Roush.
\newblock \emph{Incline algebra and applications}.
\newblock Ellis Horwood, 1984.

\bibitem[DJLC53]{dubreil53}
M.~Dubreil-Jacotin, L.~Lesieur, and R.~Croisot.
\newblock \emph{Le\c{c}ons sur la Th\'eorie des Treillis, des Structures
  Alg\'ebriques Ordonn\'ees, et des Treillis g\'eom\'etriques}, volume XXI of
  \emph{Cahiers Scientifiques}.
\newblock Gauthier Villars, Paris, 1953.

\bibitem[DK78]{dk}
S.~Dolecki and S.~Kurcyusz.
\newblock On ${\Phi}$-convexity in extremal problems.
\newblock \emph{SIAM J. Control Optim.}, 16:277--300, 1978.

\bibitem[Fan63]{fan}
K.~Fan.
\newblock On the {K}rein {M}ilman theorem.
\newblock In V.~Klee, editor, \emph{Convexity}, volume~7 of \emph{Proceedings
  of Symposia in Pure Mathematics}, pages 211--220. AMS, Providence, 1963.

\bibitem[GHK{\etalchar{+}}80]{gierzETAL}
G.~Gierz, K.~Hofmann, K.~Keimel, J.~Lawson, M.~Mislove, and D.~Scott.
\newblock \emph{A Compendium of Continuous Lattices}.
\newblock Springer, 1980.

\bibitem[GM02]{gondran02}
M.~Gondran and M.~Minoux.
\newblock \emph{Graphes, Dio\"\i des et semi-anneaux}.
\newblock TEC \& DOC, Paris, 2002.

\bibitem[Gol92]{golan92}
J.~Golan.
\newblock \emph{The theory of semirings with applications in mathematics and
  theoretical computer science}, volume~54.
\newblock Longman Sci \&\ Tech., 1992.

\bibitem[KM97]{maslovkolokoltsov95}
V.~N. Kolokoltsov and V.~P. Maslov.
\newblock \emph{Idempotent analysis and applications}.
\newblock Kluwer Acad. Publisher, 1997.

\bibitem[Kor65]{K65}
A.~A. Korbut.
\newblock Extremal spaces.
\newblock \emph{Dokl. Akad. Nauk SSSR}, 164:1229--1231, 1965.

\bibitem[LMS01]{litvinov00}
G.~Litvinov, V.~Maslov, and G.~Shpiz.
\newblock Idempotent functional analysis: an algebraic approach.
\newblock \emph{Math. Notes}, 69(5):696--729, 2001.
\newblock \doi{10.1023/A:1010266012029}.
\newblock \arxiv{math.FA/0009128}.

\bibitem[LS02]{nuclear}
G.~Litvinov and G.~Shpiz.
\newblock Nuclear semimodules and kernel theorems in idempotent analysis: an
  algebraic approach.
\newblock \emph{Doklady Math. Sci.,}, 6, 2002.
\newblock \arxiv{math.FA/0206026}.

\bibitem[Mas73]{maslov73}
V.~P. Maslov.
\newblock \emph{M\'ethodes Op\'eratorielles}.
\newblock Mir, Moscou, 1973.
\newblock Trad. fr. 1987.

\bibitem[MLRS02]{singer}
J.-E. Mart{\'{\i}}nez-Legaz, A.~M. Rubinov, and I.~Singer.
\newblock Downward sets and their separation and approximation properties.
\newblock \emph{J. Global Optim.}, 23(2):111--137, 2002.
\newblock \doi{10.1023/A:1015583411806}.

\bibitem[MLS91]{singer91}
J.-E. Mart{\'{\i}}nez-Legaz and I.~Singer.
\newblock $\vee$-dualities and $\bot$-dualities.
\newblock \emph{Optimization}, 22:483--511, 1991.

\bibitem[MS92]{maslov92}
V.~P. Maslov and S.~N. Samborski\u\i.
\newblock \emph{Idempotent analysis}, volume~13 of \emph{Advances in Soviet
  Mathematics}.
\newblock Amer. Math. Soc., Providence, 1992.

\bibitem[Rom67]{R67a}
I.~V. Romanovski{\u\i}.
\newblock Optimization of the stationary control for a discrete deterministic
  process.
\newblock \emph{Kibernetika}, 2:66--78, 1967.

\bibitem[RS01]{singer00}
A.~M. Rubinov and I.~Singer.
\newblock Topical and sub-topical functions, downward sets and abstract
  convexity.
\newblock \emph{Optimization}, 50(5-6):307--351, 2001.

\bibitem[Rub00]{rubinov}
A.~M. Rubinov.
\newblock \emph{Abstract convexity and global optimization}.
\newblock Kluwer, 2000.

\bibitem[Sin84]{singer84}
I.~Singer.
\newblock Generalized convexity, functional hulls and applications to conjugate
  duality in optimization.
\newblock In G.~Hammer and D.~Pallaschke, editors, \emph{Selected Topics in
  Operations Research and Mathematical Economics}, number 226 in Lecture Notes
  Econ. Math. Systems, pages 49--79. Springer, 1984.

\bibitem[Sin97]{ACA}
I.~Singer.
\newblock \emph{Abstract convex analysis}.
\newblock Wiley, 1997.

\bibitem[SS92]{shpiz}
S.~N. Samborski{\u\i} and G.~B. Shpiz.
\newblock Convex sets in the semimodule of bounded functions.
\newblock In \emph{Idempotent analysis}, pages 135--137. Amer. Math. Soc.,
  Providence, RI, 1992.

\bibitem[Vor67]{V67}
N.~N. Vorobyev.
\newblock Extremal algebra of positive matrices.
\newblock \emph{Elektron. Informationsverarbeit. Kybernetik}, 3:39--71, 1967.

\bibitem[Vor70]{V70}
N.~N. Vorobyev.
\newblock Extremal algebra of non-negative matrices.
\newblock \emph{Elektron. Informationsverarbeit. Kybernetik}, 6:303--311, 1970.

\bibitem[Wag91]{wagneur91}
E.~Wagneur.
\newblock Moduloids and pseudomodules. 1. {Dimension} theory.
\newblock \emph{Discrete Math.}, 98:57--73, 1991.

\bibitem[Zim76]{Zimmermann.K}
K.~Zimmermann.
\newblock \emph{Extrem\'aln\'\i\ Algebra}.
\newblock Ekonomick\'y \`ustav \u CSAV, Praha, 1976.
\newblock (in Czech).

\bibitem[Zim77]{zimmerman77}
K.~Zimmermann.
\newblock A general separation theorem in extremal algebras.
\newblock \emph{Ekonom.-Mat. Obzor}, 13(2):179--201, 1977.

\bibitem[Zim79a]{zimmer}
K.~Zimmermann.
\newblock Extremally convex functions.
\newblock \emph{Wiss. Z. P\"ad. Hochschule ``N. K. Krupskaya''}, 17:3--7, 1979.

\bibitem[Zim79b]{zimmermann}
K.~Zimmermann.
\newblock A generalization of convex functions.
\newblock \emph{Ekonom.-Mat. Obzor}, 15(2):147--158, 1979.

\bibitem[Zim81]{Zimmermann.U}
U.~Zimmermann.
\newblock \emph{Linear and Combinatorial Optimization in Ordered Algebraic
  Structures}.
\newblock North Holland, 1981.

\end{thebibliography}
\end{document}